\documentclass[aos,MSNbibl,nameyear,dvips]{arximspdf}
\usepackage{multirow,mathrsfs}
\usepackage{graphicx}

\doi{10.1214/13-AOS1195} 
\volume{42}
\issue{1}
\pubyear{2014}
\firstpage{382}
\lastpage{412}

\makeatletter
\newcommand{\rrvert}{\vert}
\newcommand{\llvert}{\vert}
\def\cal{\mathcal}

\newtheorem{proposition}{Proposition}

\newtheorem{corollary}{Corollary}
\newproclaim{assumption}{Assumption}
\newproclaim{definition}{Definition}
\newtheorem{theorem}{Theorem}

\newproclaim{example}{Example}
\newproclaim{problem}{Problem}
\def\vec{\operatorname{vec}}
\def\indep{\perp\!\!\!\!\perp}
\newcommand{\E}{\mathrm{E}}
\newcommand{\var}{\operatorname{Var}}

\newcommand{\I}{{\mathbf I}}

\newcommand{\udex}[1]{^{\nano #1}}

\newcommand{\abb}[1]{{\uppercase{#1}}}

\def\scirc{\circ}
\def\half{^{{1}/{2}}}
\def\cov{\operatorname{cov}}
\def\cid{\stackrel{\cal D}\rightarrow}
\def\real{\mathbb R}
\def\equal{=}
\newcommand{\sten}[1]{\mathscr{ #1}}

\def
\prob{{\cal P}}

\def
\prob{\mbox{\rsfsten P}}

\def\oc{^{\nano\perp}}

\newcommand{\lt}[1]{L_{\nano2}(#1)}

\def\field{{\cal F}}
\def\of{\,{ \circ}\,}
\def\T{ {T}}
\def\R{ {R}}
\def\v{\operatorname{Var}}
\def\E{ {E}}

\def\nano{\scriptscriptstyle}
\def\trans{^{\nano\mathsf T}}

\def\l2{\mbox{L}_2}

\def\inv{^{\nano-1}}
\def\T{ {T}}
\newcommand{\lo}[1]{_{\nano #1}}
\def\nano{}
\newcommand{\om}{\Omega}
\renewcommand\field[1]{{\sten B}_{\nano #1}}
\renewcommand{\real}[1]{\mathbb{R}\udex{#1}}
\newcommand{\me}[1]{\mu_{\nano #1}}
\def\tint{ \int}
\def\sgn{\operatorname{sgn}}
\newcommand{\proj}[2]{\Pi( #1 | #2 )}
\def\star{\udex*}

\newcommand{\hi}[1]{^{\nano #1}}
\def\frechet{Fr\'{e}chet }
\newcommand{\ca}[1]{{\cal #1}}
\def\om{\Omega}
\newcommand{\cts}[1]{{\ca S}_{\nano #1(Y|X)}}
\def\sdr{{SDR}}
\makeatother

\begin{document}
\begin{frontmatter}

\title{On efficient dimension reduction with respect
to a~statistical functional of interest}
\runtitle{Dimension reduction for statistical functional}
\pdftitle{On efficient dimension reduction with respect
to a statistical functional of interest}

\begin{aug}
\author[A]{\fnms{Wei} \snm{Luo}\ead[label=e1]{wzl118@psu.edu}},
\author[A]{\fnms{Bing} \snm{Li}\corref{}\thanksref{t1}\ead[label=e2]{bing@stat.psu.edu}}
\and
\author[C]{\fnms{Xiangrong} \snm{Yin}\thanksref{t2}\ead[label=e3]{xryin@stat.uga.edu}}
\thankstext{t1}{Supported in part by a NSF Grant DMS-11-06815.}
\thankstext{t2}{Supported in part by a NSF Grant DMS-12-05546.}
\runauthor{W. Luo, B. Li and X. Yin}
\affiliation{Pennsylvania State University, Pennsylvania State University
and University of Georgia}
\address[A]{W. Luo\\
B. Li\\
Department of Statistics\\
Pennsylvania State University\\
326 Thomas Building\\
University Park, Pennsylvania 16803\\
USA\\
\printead{e1}\\
\phantom{E-mail:\ }\printead*{e2}}
\address[C]{X. Yin\\
Department of Statistics\\
101 Cedar Street\\
University of Georgia\\
Athens, Georgia 30602\\
USA\\
\printead{e3}}
\end{aug}

\received{\smonth{2} \syear{2013}}
\revised{\smonth{8} \syear{2013}}

%
\begin{abstract}
We introduce a new sufficient dimension reduction framework that
targets a statistical functional of interest,
and propose an efficient estimator for the semiparametric estimation
problems of this type. The statistical functional
covers a wide range of applications, such as conditional mean,
conditional variance and conditional quantile.
We derive the general forms of the efficient score and efficient
information as well as their specific forms
for three important statistical functionals: the linear functional, the
composite linear
functional and the implicit functional. In conjunction with our
theoretical analysis, we also propose a class of one-step
Newton--Raphson estimators and show
by simulations
that they substantially outperform existing methods.
Finally, we apply the new method to construct the central mean and
central variance subspaces for
a data set involving the physical measurements and age of abalones, which
exhibits a strong pattern of heteroscedasticity.
\end{abstract}

%
\begin{keyword}[class=AMS]
\kwd{62-09}
\kwd{62G99}
\kwd{62H99}
\end{keyword}
\begin{keyword}
\kwd{Central subspace}
\kwd{conditional mean}
\kwd{variance and quantile}
\kwd{efficient information}
\kwd{efficient score}
\kwd{\frechet derivative and its representation}
\kwd{projection}
\kwd{tangent space}
\end{keyword}

\pdfkeywords{62-09, 62G99, 62H99, Central subspace, conditional mean, variance and quantile,
efficient information, efficient score, Frechet derivative and its representation, projection, tangent space}

\end{frontmatter}

\section{Introduction}\label{sec1}

The purpose of this paper is twofold: to introduce a new framework for
sufficient dimension
reduction that targets a statistical functional of interest, and to
develop semiparametrically efficient
estimators for problems of this type.

Let $X$ be a $p$-dimensional random vector and $Y$
be a random variable. In classical sufficient dimension reduction (\sdr
), we are interested in a lower dimensional linear predictor
$\zeta\trans X$, where $\zeta$ is a $p \times d$ matrix with $d < p$,
such that $Y$ is independent of $X$ given $\zeta\trans X$.
That is, the conditional distribution of $Y$ given $X$ is the same as
that of $Y$ given $\zeta\trans X$. In this problem, the identifiable
parameter is $\operatorname{span} (\zeta)$, the subspace of $\real{p}$
spanned by the columns of $\zeta$. Under mild conditions, there exists
a unique
smallest subspace that satisfies this condition, and it is called the
\emph{central subspace};
see Li (\citeyear{Li91N2,Li92}), \citet{CoWe91}, Cook
(\citeyear{Coo94,Coo96,Coo98}). For a general
discussion of the sufficient
conditions for the central subspace to exist, see \citet
{YinLiCoo08}. The \sdr\  provides us a mechanism to reduce the
dimension of the
predictor while preserving the conditional distribution of $Y$ given $X$.

\newcommand{\eff}{_{\mathrm{eff}}}
\newcommand{\neff}{_{{\nano n}, \mathrm{eff}}}
\def\opg{{OPG}}
\def\mave{{MAVE}}
\def\rmave{{RMAVE}}

In many applications, our interests are in some specific aspects of the
conditional distribution $P \lo{Y|X}$.
For example, in nonparametric regression, we are interested in the
conditional mean $E(Y|X)$; in median regression,
we are interested in the conditional median $M(Y|X)$, and in volatility
analysis we are interested in the conditional variance $\var(Y|X)$
and in supervised classification we are interested in the class label
of $Y$ given its covariates $X$.
To illuminate this point further, let us consider the model
%
\begin{equation}
Y = \mu(X \lo1) + \sigma(X \lo1 + X \lo2) \varepsilon, \label
{eq:mu sigma}
\end{equation}
where $\mu$ and $\sigma$ are unknown functions. In this case, the
central subspace, which is the two-dimensional subspace of $\real p$
spanned by $(1,0,\ldots, 0)$ and $(0, 1, 0, \ldots, 0)$,
only tells us the sufficient predictors can be any linear combination
of $X\lo1$ and~$X \lo2$,
but it does not tell us that the conditional mean is a function of $X
\lo1$, and the conditional variance is a function of $X \lo1 + X \lo
2$. Thus,
the information provided by the central subspace is clearly inadequate
if we want to build a model like (\ref{eq:mu sigma}).

Under these circumstances, it makes sense to reformulate sufficient
dimension reduction to target a specific functional, so as to provide a
more nuanced picture of the relation between $X$ and $Y$ than offered
by the central subspace. Several such efforts have been made
over the past decade or so. For example, \citet{CooLi02}
introduced the central mean subspace, which is defined by the relation
$E(Y|X) = E (Y|\alpha\trans X)$, where $\alpha$ is minimal in the same
sense as it is in the central subspace. \citet{YinCoo02}
introduced the $k$th central moment subspace through the relation
$E(Y\udex k | X) = E ( Y\udex k | \alpha\trans X)$.
\citet{ZhuZhu09} introduced the central variance subspace by
requiring that $\var(Y|X)$ is a function of $\alpha\trans X$. Zhu,
Dong and Li (\citeyear{ZhuDonLi13})
introduced a general class of estimating equations for single-index
conditional variance.
The Minimum Average Variance Estimator \mbox{(\mave)} by
\citet{Xiaetal02} also targets the central mean subspace.
\citet{KonXia12} introduced an adaptive quantile estimator for
single-index quantile
regression, which targets the conditional quantile.
It turns out the space spanned by the columns of $\alpha$ in the above
relations are all
subspaces of the central subspace. They provide refined structures for
the central subspace. As a consequence, $\alpha\trans X$ can be written
as $\beta\trans\zeta\trans X$,
and we can refine central subspace based on the sufficient predictor
$\zeta\trans X$. For convenience, we reset
$\zeta\trans X$ as $X$, and use $\tilde X$ to represent the original
predictor throughout the rest of the paper.

The first goal of this paper is to unify these problems by introducing
a general dimension reduction paradigm
with respect to statistical functional $T$
of the conditional density of $Y$ given $X$, say $\eta(x,y)$, through
the following statement:
%
\begin{equation}
\label{eq:t-sdr} T \bigl( \eta(x, \cdot)\bigr) \mbox{ is a function of $\beta\trans
x$}.
\end{equation}
Note that sufficient dimension reduction for conditional mean,
conditional moments and conditional variance discussed in
the last paragraph
are all special cases of relation (\ref{eq:t-sdr}).
The minimal subspace $\operatorname{span} (\beta)$ of $\real p$ that
satisfies this relation is called the $T$-central subspace.

The second, and the main goal of this paper is to develop
semiparametrically efficient estimators for the $T$-central subspace.
In a series of
recent papers, Ma and Zhu (\citeyear{MaZhu12,MaZhu13N1,MaZhu13N2}) use
semiparametric theory to study sufficient
dimension reduction
and develop semiparametrically efficient estimators of the central
subspace. These are related to the earlier developments by \citet{LiDon09}
and \citet{DonLi10}, which use estimating equations to relax the
elliptical distribution assumption for sufficient dimension reduction.
We extend
Ma and Zhu's approach to find semiparametrically efficient estimator
for the $T$-central subspace. We derive the general formulas
for the efficient score
and efficient information for the semiparametric family specified by
the relation (\ref{eq:t-sdr}), and further deduce their specific forms
for three important statistical functionals: the linear functionals
(L-functionals), the composite linear functionals (C-functionals)
and implicit functionals (I-functionals).
These functionals cover a wide range of applications. For example, all
conditional moments are L-functionals, all conditional cumulants
[see, e.g., \citet{McC87}] are C-functionals and quantities such as
conditional median, conditional quantile and conditional
support vector machine \cite{LiArtLi11} are I-functionals.

Using the formulas for efficient score and efficient information,
we propose a one-step Newton--Raphson
alorithm to implement
semiparametrically efficient estimation.
Compared with the semiparametric estimators of \citet{MaZhu13N1},
our algorithm has two
distinct and attractive features. First, since our algorithm relies on
the \mave-type procedure for minimization, it can be implemented by
iterations of a least squares problem without resorting to
high-dimensional search-based optimization.
Second, unlike \citet{MaZhu13N1}, our method does not require any
specific parameterization of
$\beta$ that potentially restricts the generality of their method.

The rest of the paper is organized as follows. In Section~\ref{section:cts}, we give a general formulation of sufficient
dimension reduction with respect to a statistical functional of
interest. To set the stage for further development, we lay out
the semiparametric structure of our problem in Section~\ref{section:semiparametric structure}. In Section~\ref{section:efficient score and information}, we derive the efficient
score and efficient information for
a general statistical functional.
In Sections~\ref{section:l functional}, \ref{section:c functional} and
\ref{section:i functional}, we further deduce the specific forms of the
efficient
score and efficient information for the L-, C- and I-functionals. In
Section~\ref{section:effect}, we discuss the effect of estimating the
central subspace
on the efficient score. In Section~\ref{section:estimation}, we develop
the one-step Newton--Raphson
estimation procedure
for semiparametrically efficient estimation. In Section~\ref{section:simulation}, we conduct simulation studies to compare our
method with other methods,
and in Section~\ref{section:real data} we apply our method to a data
set. Some concluding remarks are made in Section~\ref{section:discussion}. The proofs of
some technical results are given in the supplementary material [\citet{supp}].

The following notation will be consistently used throughout the rest
of the paper. The symbol $I \lo k$ denotes the $k \times k$ dimensional
identity matrix; $e \lo k$ denotes a vector
whose $k$th entry is 1 and other entries are 0;
$\indep$ indicates independence or conditional independence between two
random elements---that is, $A \indep B$ means
$A$ and $B$ are independent, and $A \indep B | C$ means $A$ and $B$ are
independent given~$C$. For integers $s$ and $t$,
$\real s$ denotes the $s$ dimensional Euclidean space, and $\real {s
\times t}$ denotes the set of $s \times t$ dimensional matrices.
For a function with multiple arguments, say $f(x,y,z)$, we use the dot
notation to represent mappings of a subset of the arguments.
For example, $f(\cdot,y,z)$ represents the mapping $x \mapsto f(x,y,z)$
where $y$ and $z$ are fixed, and $f(\cdot,\cdot, z)$ represents
the mapping $(x,y)\mapsto f(x,y,z)$ where $z$ is fixed. We use
superscripts of $X$ to index components and subscripts of $X$ to
index subjects. Thus, $X \lo i \hi j$ means the $j$th component of the
$i$th observation in a sample $X \lo1, \ldots, X \lo n$. However, $a
\hi i$ represents
power when $a$ is not $X$.

\section{Dimension reduction for conditional statistical
functional}\label{section:cts}

Let $(\om\lo X,\break  \field X, \me X)$ and $(\om\lo Y, \field Y, \me Y)$
be $\sigma$-finite measure spaces,
where $\om\lo X \subseteq\real d$ and $\om\lo Y \subseteq\real{}$ and
$\field X$
and $\field Y$ are $\sigma$-fields of Borel sets in $\om\lo X$ and
$\om\lo Y$. Let $(X,Y)$ be a pair of random elements that takes
values in $(\om\lo X \times\om\lo Y, \field X \times\field Y)$. Let
$\sten M$ be a family of densities of $(X,Y)$ with respect to $\mu=
\mu\lo X \times
\mu\lo Y$.
We assume that $\sten M$ is a \emph{semiparametric} family; that is,
there exist $\Theta\subseteq\real{r}$
and a family $\sten F$ of functions $\phi\dvtx  \om\lo X \times\om\lo Y
\times\Theta\to\real{}$ such that
\[
\sten M = \cup\{\sten M \lo\theta\dvtx \theta\in\Theta \} \qquad\mbox{where } \sten M
\lo\theta= \bigl\{ \phi(\cdot, \cdot, \theta)\dvtx \phi\in\sten F \bigr\}.
\]
Furthermore, we assume that each $\phi\in\sten F $ can be factorized
into $\lambda(x) \eta(x,y,\theta)$ where $ \lambda$ is the marginal
density of $X$,
$\eta(x, \cdot, \theta)$
is the conditional density of $Y$ given $X$. The real assumption in
this factorization is that $\theta$ appears only in the conditional density.

As an illustration, consider the single index model where
\[
Y = m \bigl(\beta\trans X\bigr) + \varepsilon,\qquad \beta\in\real d, X \indep \varepsilon,
\varepsilon\sim N\bigl(0,\sigma\hi2\bigr), X \sim N(0, \Sigma),
\]
and $m$ is an unknown function. Since $m$ is unknown, $\beta$ is
identified up to a proportional constant. To avoid the trivial case,
let us assume it has
at least one nonzero component, and further assume this is the first
component for convenience.
We can then assume without loss of generality $\beta\trans= (1,
\theta\trans)$ where $\theta\in\real{d-1}$. Then
\[
\beta= e \lo1 + \Gamma \theta,\qquad m \bigl( \beta\trans X\bigr) = m \bigl( X \hi 1 + \theta\trans
\Gamma\trans X\bigr),
\]
where $\Gamma\trans$ is the $d \times(d-1)$ matrix $(0, I \lo{d-1})$.
In this case,
$\lambda(x)$ is the p.d.f. of $N(0, \Sigma)$ and $\eta(x, \cdot,
\theta)$ is the p.d.f. of $N(m ( x \hi1 + \theta\trans\Gamma\trans
x), \sigma\hi2)$ for
a given $x$.

Now let
$
\sten F \lo1$ denote the family $\{ \eta\dvtx  \phi\in\sten F \}$, and
\[
\sten L = \{\lambda\dvtx \phi\in\sten F \},\qquad \sten H \lo\theta= \bigl\{\eta (
\cdot, \cdot, \theta)\dvtx \eta\in\sten F \lo1 \bigr\}, \qquad\sten H = \cup\{ \sten H
\lo\theta\dvtx \theta\in\Theta\}.
\]
We assume that $\sten M$ contains the true density of $(X,Y)$. That is,
there exist $\theta\lo0 \in\Theta$, $\lambda\lo0 \in\sten L$, and
$\eta\lo0 \in\sten F \lo1 $ such that
$\phi\lo0(x,y,\theta\lo0) = \lambda\lo0 (x) \eta\lo0 (x, y,
\theta\lo0)$ is
the true density of $(X,Y)$. For convenience, we abbreviate $\phi\lo0
(x,y, \theta\lo0)$, $\eta\lo0 (x, y, \theta\lo0)$, and
$\sten M \lo{\theta\lo0}$ as $\phi\lo0 (x,y)$,
$\eta\lo0 (x,y)$ and
$\sten M \lo0$.

Let $\sten G $ be a class
of densities of $Y$ with respect to $\me Y$ that contains all $\eta
(x,\cdot, \theta)$ for $\eta\in\sten F \lo1$, $\theta\in\Theta$
and $x \in\om\lo X$.
Let $T\dvtx \sten G \to\real{}$ be a mapping from $\sten G $ to $\real
{}$. Such mappings are called statistical functionals.
The functional $T$ induces the random variable
\[
x \mapsto\T\bigl(\eta(x,\cdot, \theta)\bigr),
\]
on $\om\lo X$, which we write as $T(\eta(X, \cdot, \theta))$.
Following the convention of
the conditional expectation, we write $T(\eta\lo0 ( X, \cdot, \theta
\lo0))$ as $T (Y|X)$.
This random variable can be used to
characterize a wide variety of features of a conditional density $\eta
(x, \cdot, \theta)$ that might interest us,
as detailed by the following example.

\begin{example} Let $T\dvtx \sten G \to\real{}$ be the functional $g
\mapsto\int\lo{\om\lo Y} y g(y) \,d \me Y$. Then each $\eta(\cdot,
\cdot, \theta)
\in\sten H$ uniquely
defines the mapping
\[
x \mapsto\T\bigl(\eta(x,\cdot, \theta)\bigr) = \tint\lo{\om\lo Y} y \eta (x,y,
\theta) \,d \me Y (y).
\]
That is, $\T(\eta(X, \cdot, \theta))$ is the conditional expectation
$E (Y|X)$ under $\eta(\cdot, \cdot, \theta)$. 

Let $T\dvtx \sten G \to\real{}$ be the mapping
\[
\T(g) = \tint y\udex2 g(y) \,d \me Y (y) - \biggl( \tint y g(y) \,d \me Y (y)
\biggr)\udex2.
\]
Then $T(\eta(X, \cdot, \theta))$ is the conditional variance $\var
(Y|X)$ under the conditional density $x \mapsto\eta(x, \cdot, \theta
)$. 

Let $T\dvtx \sten G \to\real{}$ be the functional defined by the equation
in $m$
%
\begin{equation}
\label{eq:median} \tint\sgn(y - m) g (y) \,d \me Y(y) = 0,
\end{equation}
where $\sgn(a)$ is the sign function that takes the value 1 if $a \ge
0$ and takes the value $-1$ if $a < 0$.
The solution to (\ref{eq:median}) is the median of $Y$. Each $\eta(\cdot
, \cdot, \theta) \in\sten H$ uniquely
defines the mapping $T ( \eta(X, \cdot, \theta))$, which is the
conditional median of $Y$ given $X$ under the conditional density
$x \mapsto\eta(x, \cdot, \theta)$.
\end{example}

We now give a rigorous definition of the $T$-central subspace.

\begin{definition}\label{definition:sdr} If there is a matrix $\gamma
\in\real{d \times u}$,
with $u < d$, such that
$
T (\eta\lo0 (X, \cdot, \theta\lo0))$ is measurable with respect to
$\sigma(\gamma\trans X)$,
then we call $\operatorname{span} (\gamma)$ a sufficient dimension
reduction subspace for $T$. The intersection of all such spaces is
called the central
subspace for conditional functional $T$, or the $T$-central subspace.
\end{definition}

\def\cs{{\ca S}\lo{Y|X}}

We denote the $T$-central subspace by $\cts T$.
For example, if $T$ is the conditional mean functional, then $\cts T$
becomes the central mean subspace, which we write as $\cts E$;
if $T$ is the conditional median functional, then $\cts T$ becomes the
central median subspace, which we write as $\cts M$;
if $T$ is the conditional variance functional, then $\cts T$ becomes
the central variance subspace, which we write as $\cts V$.
It is easy to see that $\cts T \subseteq\cs$: this is because $Y
\indep X | \beta\trans X$ implies
\[
T(Y|X) = E \bigl[ T (Y|X) | X \bigr] = E \bigl[ T(Y|X) | X, \beta\trans X\bigr]
= E \bigl[ T(Y|X) |\beta\trans X\bigr].
\]
In the following, we denote the dimension of the $\cts T$ by $s$ and
any basis matrix of $\cts T$ (of dimension $d \times s$) as $\beta$.

The semiparametrically efficient score derived by \citet
{MaZhu13N1} corresponds to the case where $T$ exerts no restriction on
$\eta$; that is
$T (\eta) = \eta$. However, technically this is not covered by
Definition~\ref{definition:sdr} as a special case, because the range of
$T$ considered in this paper is always a Euclidean space.

\section{Formulation of the semiparametric problem}\label
{section:semiparametric structure}
\newcommand{\effsco}[1]{S_{\mathrm{ eff}}(#1)}
\newcommand{\effscot}[1]{S_{\mathrm{ eff}}\trans(#1)}
\newcommand{\effinfo}[1]{J_{\mathrm{ eff}}(#1)}
\newcommand{\effinfoi}[1]{J_{\mathrm{ eff}}\inv(#1)}
\newcommand{\putop}[1]{\stackrel{\nano #1}}
\def\ka{\kappa}
\newcommand{\at}[2]{|\lo{#1 = #2}}

To set the stage for further development, we first outline
the basic semiparametric structure of our problem.
Let
$
\lt{\phi\lo0}=\{ r\dvtx \tint r\udex2 \phi\lo0 \,d \mu< \infty\}$. Let
$\langle\cdot, \cdot\rangle\lo{\phi\lo0}$ and $\| \cdot\| \lo
{\phi\lo0}$ denote the inner product and norm in $\lt{\phi\lo0}$.
For a technical reason, it is easier to work with an embedding of
$\sten M$ into $\lt{\phi\lo0}$, defined as
\[
R\dvtx \phi\mapsto 2 \bigl(\phi\half- \phi\lo0 \half\bigr)/\phi\lo0 \half \equiv r.
\]
Let $\R(\sten M)=\{\R(\phi)\dvtx \phi\in\sten M\}$. This
transformation ensures that $R (\sten M)\subseteq\lt{\phi\lo0}$; whereas
additional assumptions are needed to ensure $\sten M \subseteq\lt
{\phi\lo0}$.
Also note that $\R(\phi\lo0)$ is the 0 element in $\lt{\phi\lo0}$.

A \emph{curve} in $\R(\sten M \lo0)$ that passes through $r\lo0 = 0$
is any mapping $\alpha\mapsto r\lo\alpha(\cdot)$ from $[0, 1) \to
R(\sten M \lo0)$
that is \frechet
differentiable at $\alpha= 0$. That is, there is a member $\dot r \lo
0$ of $\lt{\phi\lo0}$ such that
\[
\| r \lo\alpha- r \lo0 - \dot r \lo0 \alpha\|\lo{ {\phi\lo0}} = o\bigl (|\alpha|\bigr).
\]
The \emph{tangent space} $\sten T \lo\phi$ of $\R(\sten M \lo0)$ at
$r\lo0$
is the closure of the subspace of $\lt{\phi\lo0}$ spanned by $\dot
r\lo0$ along all curves.

Let ${\putop\circ r} \lo0 \in[\lt{\phi\lo0}]\udex r$ be the score
with respect to $\theta$; that is,
\[
\bigl\| R \bigl( \phi\lo0 ( \cdot, \theta) \bigr)
- R \bigl( \phi\lo0 (\cdot, \theta
\lo0)\bigr) -{ \putop\circ r \lo0} {}^{\nano{\mathsf T}} (\theta- \theta \lo0) \bigr\|
\lo {
{\phi\lo0}} = o\bigl ( \| \theta- \theta\lo0 \| \bigr).
\]
Let $\proj{{\putop\circ r} \lo0} {\sten T \lo\phi\oc}$ be the
componentwise projection of the random vector ${\putop\circ r} \lo0$
on to the orthogonal
complement of the
tangent space $\sten T \lo\phi$. This projection is called the
efficient score, and we denote it by $\effsco{X,Y, \theta\lo0}$. The matrix
\[
\effinfo{\theta\lo0} = E \bigl[ \effsco{X,Y, \theta\lo0} \effscot {X,Y, \theta\lo0} \bigr]
\]
is called the efficient information. Now let $(X \lo1, Y \lo1), \ldots
, (X \lo n, Y \lo n)$ be an i.i.d. sample of $(X,Y)$. For a function
$h$ of $(x,y)$,
let $E \lo n h(X,Y)$ denote the sample average of $h(X \lo1, Y \lo1),
\ldots, h(X \lo n, Y \lo n)$.
Under some conditions, if $\hat\theta$ is the solution to the
estimating equation
%
\begin{equation}
\label{eq:score equation} E \lo n \effsco{X,Y, \theta} = 0,
\end{equation}
then
$
\sqrt n ( \hat\theta- \theta\lo0 ) \cid N(0, \effinfoi{\theta\lo
0} )$.
Moreover, for any estimate $\tilde\theta$ of $\theta\lo0$ that is
regular with respect to $\sten T \lo\phi$,
$
\sqrt n ( \tilde\theta- \theta\lo0 )$ can be decomposed as $\sqrt n
(\hat\theta- \theta\lo0) + \Delta\lo n
$
where two terms are asymptotically independent. This result, well known
in the semiparametric literature as the H\'{a}jek--LeCam
convolution theorem, implies $\hat\theta$
has the smallest asymptotic variance among all regular estimators with
respect to $\sten T \lo\phi$. That is, $\hat\theta$ is
semiparametrically efficient.
For a comprehensive exposition of this theory,
see Bickel et al. [(\citeyear{Bicetal93}), Chapter~3] or van der Vaart
[(\citeyear{van98}), Chapter~25].

We now investigate how the sufficient dimension reduction in Definition~\ref{definition:sdr} specifies the semiparametric family $\sten M$,
and what is the meaning of the parameter $\theta$ in this context.
Since our goal is to estimate $\operatorname{span} (\beta)$, we need
fewer parameters than $ds$.
In fact, the set $\{\operatorname{span} (\beta)\dvtx \beta\in\real{d
\times s}, \mbox{rank}(\beta) = s\}$ is a Grassmann manifold, which has
dimension $s(d-s)$ [see, e.g., Edelman, Arias, and Smith (\citeyear
{EdeAriSmi99})].
There always exists a smooth parameterization $\beta= \beta(\theta)$, where
$\theta\in\real{s (d-s)}$, because $\operatorname{span} (\beta)$ is
determined if a certain $s \times s$ submatrix of $\beta$ is
fixed as $I \lo s$ and the complementary $(d-s)\times s$ block has free
varying entries.
The specific form of the parameterization is not important to us.

Let $\sigma\lo\theta(X)$ be the $\sigma$-field generated by $\beta
\trans(\theta) X$.
Because $T ( \eta(X, \cdot, \theta))$ is measurable with respect to
$\sigma\lo\theta(X)$ if and only if
\[
T \bigl( \eta(X, \cdot, \theta)\bigr) = \E\bigl[ T \bigl( \eta(X, \cdot, \theta)
\bigr)|\sigma \lo\theta(X)\bigr],
\]
the semiparametric family for our purpose is $\sten M = \cup\{ \sten M
\lo\theta: \theta\in\real{s (d-s)}\}$, where
\[
\sten M \lo\theta= \bigl\{\phi(\cdot, \cdot, \theta)\dvtx \phi\in\sten F, \T
\bigl(\eta (x,\cdot, \theta)\bigr) = \E\lo\lambda\bigl[ \T\bigl( \eta (X,\cdot,
\theta ) \bigr) | \sigma\lo\theta(X) \bigr]\lo x \mbox{ $\forall x \in \om\lo X$}
\bigr\}.
\]
Here, for a sub-$\sigma$-field $\sten A$ of $\field X$ and a function
$f(X,Y)$, $E[f(X,Y)|\sten A] \lo x$ denotes the evaluation of
the conditional expectation $E[f(X,Y)|\sten A]$ at $x$.

In this paper, we will focus on the development of the efficient score,
the efficient information and an accompanying estimation procedure, but
will not give a rigorous proof of the asymptotic results [including the
asymptotic distribution of $\sqrt n (\hat\theta- \theta\lo0)$
and the convolution theorem], because it would far exceed the scope of
the paper and
because we do not expect the proof will fundamentally deviate from that
given in \citet{MaZhu13N1}. In addition, as mentioned earlier,
by design, our method is applied to the sufficient predictor
corresponding to the central subspace, whose dimension is relatively low.
Consequently, we expect no surprises as regards the validity of $\sqrt
n$-rate of convergence
of our estimator. In the meantime, our simulation studies provide
strong evidence that the our efficient estimator does approach
the theoretical semiparametric variance bound for modestly large sample
sizes.\looseness=1

\section{Efficient score and efficient information}\label
{section:efficient score and information}

In this section, we derive the efficient score and efficient
information for the semiparametric problem set up
in Section~\ref{section:semiparametric structure}.
To this end, we first derive
the tangent space $\sten T \lo\phi$ for a fixed $\theta\in\Theta$.
Let $\sten T \lo\eta$ be the tangent space of $R(\sten H \lo\theta)$
at $R(\eta\lo0 (\cdot, \cdot, \theta)) = 0$, and
$\sten T \lo\lambda$
be the tangent space of $R(\sten L)$ at $R(\lambda\lo0) = 0$.

\newcommand{\tangent}[1]{\sten T \lo #1}

\begin{proposition}\label{lemma:sum} The following relations hold:
\begin{longlist}
\item[1.]$\sten T \lo\phi= \sten T \lo\eta+ \sten T \lo\lambda$;

\item[2.]$\sten T \lo\eta\perp\sten T \lo\lambda$ in terms of the inner
product in $\lt{ \phi\lo0 }$;

\item[3.]$\sten T \lo\phi\oc= \sten T \lo\lambda\oc\cap\sten T \lo
\eta\oc$.
\end{longlist}
\end{proposition}

This proposition was verified and used in Ma and Zhu (\citeyear
{MaZhu12}, \citeyear{MaZhu13N1}).
Since the family $\sten L$ has no constraint, its tangent space is
straightforward, as given in the next proposition, which is
taken from Bickel et al. [(\citeyear{Bicetal93}), page 52].

\begin{proposition}\label{proposition:tangent lambda} $\tangent\lambda
$ consists of all functions $h$ in $\lt{\lambda\lo0}$ with\break
$E\lo{\lambda\lo0} h(X) = 0$.
\end{proposition}

To compute $\tangent\eta$, we introduce a new functional for each
fixed $x \in\om\lo X$. Let $\sten H \lo{\theta,x}$ be the class of densities
$\{\eta(x, \cdot, \theta)\dvtx \eta\in\sten H \lo\theta\}$.
Let $R \lo x$ denote the mapping
%
\begin{eqnarray}
\label{eq:r embedding} &&R \lo x\dvtx  \sten H \lo{\theta, x} \to
L_{\nano2}\bigl(\eta\lo0 (x, \cdot, \theta ) \bigr),
\nonumber
\\[-8pt]
\\[-8pt]
\nonumber
&&\qquad\eta(x, \cdot, \theta) \mapsto2 \bigl[\eta\half(x, \cdot) - \eta\lo0\half(x,
\cdot, \theta) \bigr] / \eta\lo0 \half(x, \cdot, \theta).
\end{eqnarray}
This mapping is different from $R$, which is from $\sten M$ to $\lt
{\phi\lo0}$. Nevertheless, note that $R (\eta(\cdot, \cdot, \theta)) (x,y)
= R \lo x (\eta(x,\cdot, \theta))(y)$.
Let
%
\begin{equation}
\label{eq:tx} \T\lo x\dvtx
L_{\nano2}\bigl(\eta\lo0 (x, \cdot, \theta)\bigr)
 \to\real{},\qquad r(x,
\cdot, \theta) \mapsto\T\of R\lo x\inv\bigl( r(x,\cdot, \theta)\bigr).
\end{equation}
The \frechet derivative of $T \lo x$ at $r(x, \cdot, \theta)$ is
denoted by $\dot T \lo x ( r (x, \cdot, \theta))$. This is a bounded linear
functional from $\lt { \eta\lo0 (x, \cdot, \theta) }$ to $\real{}$.

\begin{theorem}\label{theorem:tangent eta} Suppose, for each $x \in\om
\lo X$, the functional $\T\lo x$ is \frechet differentiable at 0.
Let $\tau(x, \cdot, \theta)$
be the
Riesz representation of $\dot T \lo x (0)$ and assume $\tau (\cdot,
\cdot, \theta) \in\lt{\phi\lo0}$. Then
%
\begin{eqnarray}\label{eq:or tan}
\tangent\eta&\subseteq&\bigl\{ \bigl[ h(x) - \E\bigl(h(X) |\sigma\lo\theta(X)
\bigr)\lo x \bigr] \tau(x,y, \theta) + g(x)\dvtx h, g \in \lt{\lambda\lo0}
\bigr\}
\oc
\nonumber
\\[-8pt]
\\[-8pt]
\nonumber
&\equiv& \sten U\oc.
\end{eqnarray}

Moreover, if, for each $ x \in\om\lo X$,
the function $r(x, \cdot, \theta) \mapsto T \lo x (r(x, \cdot, \theta
))$ is continuously \frechet differentiable in a neighborhood of
$0 \in\lt{\eta\lo0 (x, \cdot, \theta)}$, then
$
\tangent\eta\supseteq\sten U\oc$.
\end{theorem}

\def\Re{\mathbb{R}}
\mathversion{normal}

The proof of this theorem is technical and is presented in the supplementary material [\citet{supp}] (Section \abb{i}).
From Theorem~\ref{theorem:tangent eta} and Propositions \ref{lemma:sum}
and \ref{proposition:tangent lambda},
we can easily derive the form of $\tangent\phi\oc$, as follows.
%
\begin{corollary} Under the assumptions of Theorem~\ref{theorem:tangent eta},
\[
\sten T \lo\phi\oc= \bigl\{ \bigl[ h(x) - E \bigl( h(X) | \sigma\lo\theta(X)
\bigr)\lo x \bigr] \bigl[ \tau(x, y, \theta) - E\bigl( \tau(X,Y, \theta)|X\bigr)
\lo x \bigr]\dvtx h \in \lt{\lambda\lo0} \bigr\}.
\]
\end{corollary}

\newcommand{\score}[2]{\stackrel{\nano\circ}{r}\lo0 (#1, #2, \theta\lo0)}

We now compute the efficient score, which is the projection of the true
score with respect to $\theta$ on to $\sten T \lo\phi\oc$. Let
\[
r \lo0 (x, y, \theta) = 2 \bigl[ \phi\lo0 \udex{1/2} (x, y, \theta ) - \phi\lo0
\udex{1/2}(x,y, \theta\lo0 ) \bigr] / \phi\lo0 \udex{1/2} (x,y, \theta\lo0).
\]
The true score for the parameter of interest is the \frechet derivative
\[
\partial r \lo0 ( x, y, \theta) / \partial\theta \at\theta{\theta \lo0}.
\]
To differentiate from $\dot r \lo0 (x, y, \theta\lo0)$, we denote
the above derivative by
$\stackrel{\nano\circ}{r} \lo0 (x, y, \theta\lo0)$. This is an $s
(d-s)$-dimensional vector.
Since the mapping $T \lo x\dvtx \lt{ \eta\lo0 (x, \cdot, \theta) } \to
\real{}$ and $R \lo x$ also depend on $\theta$, we now write them
as $T \lo{x, \theta}$ and
$R \lo{x, \theta}$. We use $T \lo x$ to denote
the mapping $T \lo{x, \theta\lo0}$.
Following Bickel et al. [(\citeyear{Bicetal93}), Chapter~3], we use
$\proj{f} {\sten A}$ to represent the projection of a function $f$ on
to a subspace $\sten A$ of
$\lt{\phi\lo0}$.

\begin{theorem}\label{theorem:general e score} Suppose the following
conditions hold:
\begin{longlist}[1.]
\item[1.] For each $x \in\om\lo X$ and $\theta$ in a neighborhood of
$\theta\lo0$, the mapping
\[
T \lo{x, \theta} \dvtx
L_{\nano2}\bigl(\eta\lo0 ( x, \cdot, \theta )\bigr)
 \to\real{}
\]
is continuously \frechet differentiable in a neighborhood
of $0 \in\lt{\eta\lo0 (x, \cdot, \theta ) }$. Let $\tau(x, y,
\theta)$ be the Riesz representation of $\dot T \lo{x, \theta} (0)$ and
$\tau\lo c (x,y,\theta)$ be its centered version
$
\tau(x,y,\theta) - \E\lo\theta[ \tau(X, Y, \theta ) | X ]\lo x$.
\item[2.] The function $\theta\mapsto r \lo0 (x, \cdot, \theta)$ is
\frechet differentiable at $\theta= \theta\lo0$ with \frechet derivative
$\putop\scirc r \lo0 (x, \cdot, \theta\lo0)$.
\item[3.] If
\begin{eqnarray*}
q \lo1 (x, \theta\lo0) &\equal& \E \bigl[\score X Y \tau\lo c (X, Y, \theta\lo0) |
X\bigr] \lo x,
\\
q \lo2 (x, \theta\lo0) &\equal& \E\bigl[\tau\udex2 \lo c(X, Y, \theta\lo0) | X
\bigr] \lo x,
\\
q \lo3 (x, \theta\lo0) &\equal& q \lo1 (x,\theta\lo0) \\
&&{}- \E\lo {\theta\lo0} \bigl[q
\lo1 (X,\theta\lo0) q \lo2 \inv (X,\theta\lo0) | \sigma\lo{\theta\lo0} (X)
\bigr]_x/\E\bigl[q \lo2 \inv (X, \theta\lo0) | \sigma\lo{\theta\lo0}
(X)\bigr]_x,
\\
q \lo4 (x, \theta\lo0 ) &\equal &q \lo2\inv (x,\theta\lo0 ) q \lo 3(x, \theta
\lo0),
\end{eqnarray*}
where $q \lo2 \inv(x, \theta\lo0)$ is the reciprocal of $q \lo2
(x, \theta\lo0 )$, then $q \lo4 (x, \theta\lo0) \in\lt{ \lambda
\lo0}$.
\end{longlist}
Then
%
\begin{equation}
\effsco{x, y, \theta\lo0} =
\Pi\bigl( \score x y | \tangent\phi\oc \bigr)
 = q
\lo4(x, \theta\lo0) \tau\lo c (x,y,\theta\lo0). \label
{eq:projected score} %
\end{equation}
\end{theorem}

\def\star{^{\nano*}}

\begin{pf} Let $u \udex* (x, y, \theta\lo0) = q \lo4(x, \theta\lo
0) \tau\lo c (x,y,\theta\lo0)$.
By the projection theorem, it suffices to show:
\begin{longlist}[(a)]
\item[(a)] $u \udex* (\cdot, \cdot, \theta\lo0) \in\tangent\phi\oc$;
\item[(b)] for any $u \in\tangent\phi\oc$,
%
\begin{equation}
\bigl\langle\score\cdot\cdot, u \bigr\rangle\lo{\phi\lo0} = \bigl\langle u\star(\cdot,
\cdot, \theta\lo0 ), u \bigr\rangle\lo{\phi\lo0}. \label
{eq:projection}
\end{equation}
\end{longlist}
By condition 3, $q \lo4 (\cdot, \theta\lo0) \in\lt{\lambda\lo0}$.
Moreover, by the definition of $q \lo4$ in condition~3 it is easy to
verify that
$\E (q \lo4 (X, \theta) |\sigma\lo{\theta\lo0} (X) ) = 0$. Hence,
assertion (a) holds.

Because $u\in\tangent\phi\oc$, it has the form
$
h(x, \theta\lo0) \tau\lo c (x, y, \theta\lo0)
$ for some $h(\cdot, \theta\lo0) \in\lt{\lambda\lo0}$ satisfying
$\E[ h(X, \theta\lo0) | \sigma\lo{\theta\lo0} (X)] = 0$. Hence,
the right-hand side of (\ref{eq:projection}) is
\begin{eqnarray*}
\E\lo\theta \bigl[h(X, \theta\lo0) q \lo4 (X, \theta\lo0) \tau\lo c \udex2 (X,
Y, \theta\lo0) \bigr] &\equal& \E \bigl[h(X, \theta\lo0) q \lo4 (X, \theta\lo0 ) q
\lo2 (x, \theta\lo0 )\bigr]
\\
&\equal& \E \bigl[h(X, \theta\lo0) q \lo3 (X, \theta\lo0 ) \bigr].
\end{eqnarray*}
Substitute the definition of $q \lo3 (x, \theta\lo0 )$ into the
right-hand side, and it becomes
\[
\E\lo\theta \bigl\{h(X, \theta) \bigl\{q \lo1 (x,\theta\lo0) - \E\bigl[q \lo1 (x,
\theta\lo0 ) q \lo2 \inv (x,\theta\lo0 )|\sigma\lo{\theta \lo0} (X)\bigr]/\E
\bigl[q \lo2 \inv(x,\theta\lo0 ) |\sigma\lo{\theta\lo0 } (X)\bigr] \bigr\}\bigr\}.
\]
However, because $\E( h (X, \theta\lo0) | \sigma\lo{\theta\lo
0}(X)) = 0$, the equation above reduces to
\[
\E \bigl[h(X, \theta\lo0 ) q \lo1 (X,\theta\lo0) \bigr]= \E \bigl[h(X, \theta
\lo0) \putop\circ r \lo0 (X, Y, \theta\lo0 ) \tau \lo c (X, Y, \theta\lo0 )
\bigr],
\]
which is the left-hand side of (\ref{eq:projection}).
\end{pf}

The next corollary, which follows directly from Theorem~\ref
{theorem:general e score}, gives the general form for the efficient information
estimating $\cts T$.

\begin{corollary} Under the assumptions of Theorem~\ref{theorem:general
e score}, the efficient information for estimating $\cts T$
is given by
%
\begin{equation}
\label{eq:effinfo} \effinfo{\theta\lo0} = \E \bigl[ q \lo3 (X, \theta\lo0 ) q \lo3
\trans (X, \theta\lo0 ) q \lo2 \inv (X,\theta\lo0 ) \bigr].
\end{equation}
\end{corollary}

In the next three sections, we apply the general result in Theorem~\ref
{theorem:general e score} to derive the explicit
forms of the
efficient scores for three types of commonly used statistical
functionals: the linear functionals, the composite linear functionals
and the implicit functionals.
The common thread that runs through these developments is the
calculation of the Riesz representation $\tau(x,y, \theta\lo0)$ of
the \frechet derivative
$\dot T \lo x (0)$.

\section{Linear statistical functionals}\label{section:l functional}

\def\gateaux{{G\^{a}teaux }}

Dimension reduction of this type is the direct generalizations of the
central mean subspace [\citet{CooLi02}] and the central
moment subspace [\citet{YinCoo02}]. It can also be viewed as a
generalization of the single- and multiple-index models
[see, e.g., H\"{a}rdle, Hall and Ichimura (\citeyear
{HarHalIch93})]. Let $f \dvtx \om\lo Y \to\real{}$ be a square-integrable
function. Let $L$ be the functional
\[
L\dvtx \sten G \to\real{}, \qquad g \mapsto\tint\lo{\om\lo Y} f(y) g (y) \,d \me Y(y).
\]
The corresponding conditional statistical functional is
\begin{eqnarray*}
&& L \lo{x, \theta} \bigl(r (x, \cdot, \theta )\bigr)\\
&&\qquad \equiv L \of R \lo{x, \theta
}\udex{-1} \bigl(r (x, \cdot, \theta )\bigr) = \tint\lo{\om\lo Y} f(y) \bigl(1 +
r(x,y, \theta )/2\bigr)\udex2 \eta\lo0 (x, y, \theta ) \,d y.
\end{eqnarray*}
The $L$-central subspace is defined by the relation
%
\begin{equation}
\label{eq:l central} \E\bigl[ f(Y) | X \bigr] = \E\bigl[ f(Y) |\sigma\lo{\theta\lo0}
(X)\bigr].
\end{equation}

\begin{theorem} Suppose the conditions 1, 2, 3 in Theorem~\ref
{theorem:general e score} are satisfied for $L \lo{x, \theta}$.
Then the efficient score for $\theta$ in problem (\ref{eq:l central})
is given by (\ref{eq:projected score}) in which
\begin{eqnarray*}
\tau\lo c(x,y,\theta\lo0) &\equal& f(y) - E \lo{\theta\lo0} \bigl[ f(Y) | \sigma
\lo{\theta\lo0} (X) \bigr],
\\
q \lo1 (x,\theta\lo0 )& \equal& {\partial\E\lo\theta \bigl[ f(Y) | \sigma\lo\theta
(X)\bigr] \lo x}/{ \partial\theta} \at\theta{ \theta \lo0},
\\
q \lo2 (x,\theta\lo0 ) &\equal& {\E\lo{\theta\lo0}\bigl[ f\udex2 (Y) | X \bigr]\lo
x - \E\lo{\theta\lo0}\udex2 \bigl[f (Y) | \sigma\lo{\theta \lo0} (X)\bigr] \lo
x}.
\end{eqnarray*}
\end{theorem}

\begin{pf}
Because $L \lo x$ is \frechet differentiable at $0$, its \frechet
derivative is the same as the \gateaux derivative [\citet{Bicetal93},
page 455], which
is defined by
\[
r (x, \cdot, \theta\lo0 ) \mapsto\partial L \lo x \bigl( \epsilon r (x, \cdot,
\theta\lo0 )\bigr) / \partial\epsilon\at\epsilon0.
\]
However, because
\begin{eqnarray*}
&&\partial L \lo x \bigl( \epsilon r (x, \cdot, \theta\lo0 )\bigr) / \partial
\epsilon\at\epsilon0\\
&&\qquad\equal \partial \biggl[ \tint\lo{\om\lo Y} f(y) \bigl(1 +
\epsilon r(x,y, \theta\lo 0 )/2\bigr)\udex2 \eta\lo0 (x, y, \theta\lo0) \,d y
\biggr]\Big/ \partial\epsilon \at\epsilon0
\\
&&\qquad\equal\tint\lo{\om\lo Y} f(y) r (x,y, \theta\lo0 ) \eta\lo0 (x, y, \theta\lo0)
\,d y = \bigl\langle f, r(x, \cdot, \theta\lo0) \bigr\rangle\lo{\eta \lo0 (x,
\cdot, \theta\lo0)},
\end{eqnarray*}
the Riesz representation of $\dot T \lo x (0)$ is $f$. Hence, by
Theorem~\ref{theorem:general e score},
\[
q \lo2 (x, \theta\lo0 ) = { \E\bigl[\tau\udex2 \lo c (X, Y, \theta\lo 0 ) | X
\bigr] \lo x} = {\E\bigl[ f\udex2 (Y) | X \bigr]\lo x - \E\udex2 \bigl[f (Y) |
\sigma\lo {\theta\lo0} (X)\bigr]\lo x}.
\]
%
Also, for each $\theta$,
\begin{eqnarray*}
&&\tint f(y) \bigl(1 + r \lo0 (x, y, \theta)/2\bigr)\udex2 \eta\lo0 (x, y, \theta
\lo0) \,d \me Y(y)\\
 &&\qquad\equal \tint f(y) \eta\lo0 ( x, y, \theta) \,d \me Y (y)
\\
&&\qquad\equal\E\lo\theta\bigl[f(Y)|X\bigr]\lo x = E \lo\theta\bigl[ f(Y) | \sigma\lo
\theta(X) \bigr] \lo x.
\end{eqnarray*}
Take \frechet derivative with respect to $\theta$ on both sides to obtain
\[
q \lo1 (x,\theta\lo0 ) = \E\bigl[ f(Y) \putop\circ r \lo0 (X, Y, \theta\lo0 ) |
X\bigr]\lo x = \partial\E\lo\theta\bigl[ f(Y) | \sigma\lo \theta(X)\bigr] \lo x /
\partial\theta \at\theta{ \theta\lo0}
\]
as desired.
\end{pf}


\begin{example} The central mean subspace introduced by \citet
{CooLi02} is a special case of the $L$-central subspace
with $f(y) = y$. The efficient score and efficient information are
given by (\ref{eq:projected score}) and (\ref{eq:effinfo}) where
%
\begin{eqnarray}
\label{eq:3 formulas} %
\tau\lo c (x, y, \theta\lo0 ) &\equal& y - \E\lo{\theta
\lo0} \bigl[Y | \sigma\lo{\theta\lo0} (X)\bigr]\lo x,
\nonumber\\
q \lo1 (x,\theta\lo0 )& \equal&\partial\E\lo\theta \bigl[ Y | \sigma \lo{\theta }
(X)\bigr] \lo x / \partial\theta\at\theta {\theta\lo0},
\\
q \lo2 (x,\theta\lo0 ) &\equal&\v\lo{\theta\lo0} (Y | X)_x = \E \lo{
\theta\lo0}\bigl(Y\udex2 | X\bigr) - \E\lo{\theta\lo0}\udex2 \bigl[Y | \sigma\lo{\theta\lo0}
(X)\bigr].\nonumber %
\end{eqnarray}
For example, if the central mean subspace has dimension 1 and is
spanned by $c + \Gamma\theta\lo0$ for some $c \in\real {p}$ and
$\Gamma\in\real{p \times(p-1)}$, as described in the second
paragraph of Section~\ref{section:cts}, then
\begin{eqnarray*}
\tau\lo c (x, y, \theta\lo0 ) &\equal& y - \E\lo{\theta\lo0}\bigl (Y |X \hi1 + \theta
\lo0 \trans\Gamma\trans X \bigr)\lo x,
\\
q \lo1 (x,\theta\lo0 ) &\equal&\Gamma\trans x \bigl[ \partial\E\lo \theta\bigl ( Y |X
\hi1 + \theta\trans\Gamma\trans X \bigr) \lo x / \partial \bigl( \theta\trans\Gamma\trans x\bigr)
\at\theta {\theta\lo0} \bigr].
\end{eqnarray*}
Therefore, the efficient score is
\begin{eqnarray*}
S\eff (x, y, \theta)& = & \frac{1}{ \var\lo \theta (Y | X)\lo x} \frac{\partial E \lo{\theta} (Y |X \hi1 + \theta\trans\Gamma\trans
X)\lo x} {\partial(\theta\trans\Gamma x)}
\\
&&{}\times \biggl\{x - \frac{E \lo\theta[X / \var\lo\theta(Y | X) |X \hi1
+ \theta\trans\Gamma\trans X]\lo x } {
E \lo\theta [1 / \var\lo\theta (Y | X) |X \hi1 + \theta\trans
\Gamma\trans X]\lo x} \biggr\} \\
&&{}\times \bigl[y - E\lo{\theta} \bigl(Y |X \hi1 +
\theta \trans\Gamma\trans X\bigr)\lo x\bigr].
\end{eqnarray*}
The efficient information is
\begin{eqnarray*}
J\eff(\theta)& =& E \lo\theta\biggl[\frac{1}{\var\lo\theta( Y
| X )} \biggl(
\frac{\partial E\lo{\theta} (Y |X \hi1 + \theta\trans\Gamma
\trans X) }{\partial(\theta\trans\Gamma\trans X )} \biggr) \hi2
\\
&&\hspace*{2pt}\quad{}\times \biggl(X - \frac{E \lo\theta [X / \var(Y | X) |X \hi1 + \theta
\trans\Gamma\trans X] }{ E \lo\theta [1 / \var\lo\theta(Y | X) |X
\hi1 + \theta\trans \Gamma\trans X]} \biggr) \hi{\otimes2}\biggr],
\end{eqnarray*}
where $A\udex{\otimes2}$ denotes $A A\trans$ for a matrix $A$. 

Alternatively, the efficient score and information can be written in
the original parameterization $\beta$.
See the supplementary material [\citet{supp}] (Section \abb{ii}) for their explicit forms
in the $\beta$-parameterization.

The central $k$th moment space [\citet{YinCoo02}] is a special case of the
L-functional with $f(y) = y\udex k$. The efficient score
and efficient information where~(\ref{eq:projected score}) $q \lo1$,
$q \lo2$ and $\tau\lo c$ given
by formulas similar to (\ref{eq:3 formulas}) with $Y$ replaced by
$Y\udex k$.

\citet{ZhuZen06} considered a dimension reduction problem defined
through the characteristic function
$
\E\lo{\theta\lo0}( e\udex{i t Y } | X ) = \E\lo{\theta\lo0}[
e\udex{i t Y } | \sigma\lo{\theta\lo0}(X) ]$.
They used this relation to recover the central space, but if our goal
is to estimate $\theta$ defined through
this relation, then
\begin{eqnarray*}
\tau\lo c (x, y, \theta\lo0) &\equal& e \udex{i t y} - \E\lo{\theta \lo0}\bigl[ e
\udex{i t Y}| \sigma\lo{\theta\lo0} (X)\bigr]\lo x,
\\
q \lo2 (x,\theta\lo0) &\equal&\v\lo{\theta\lo0}\bigl(e\udex{i t Y} | X\bigr)\lo x,
\\
q \lo1 (x,\theta\lo0 ) &\equal& \partial\E\lo\theta \bigl[ e\udex{i t Y} | \sigma
\lo\theta (X)\bigr] \lo x / \partial\theta \at\theta {\theta\lo0}.
\end{eqnarray*}
The efficient score can be obtained by substituting the above into (\ref
{eq:projected score}).
\end{example}

\section{Composite linear statistical functionals}\label{section:c
functional}

\def\inv{\udex{-1}}
\newcommand{\tsum}[3]{ \sum_{\nano #1 = #2}\udex{#3}}

We now consider a nonlinear function of several linear functionals,
which is motivated by dimension reduction for
conditional variance considered in \citet{ZhuZhu09} and the
single-index conditional heteroscedasticity model in Zhu, Dong and Li
(\citeyear{ZhuDonLi13}).
See also \citet{XiaTonLi02}.
In fact, all cumulants are functionals of this type.
Let $T\lo1, \ldots, T\lo k$ be bounded linear functionals from $\sten
G$ to $\real{}$. That is,
\[
\T\lo\ell(g) = \tint f \lo\ell(y) g(y) \,d \me Y(y),\qquad \ell= 1, \ldots, k,
\]
where $f_1, \ldots, f_k$ are square-integrable with respect to any
density $g\in\sten G$.
Let $\rho: \real{k} \to\real{}$ be a differentiable function.
Then
$
C\dvtx g \mapsto\rho( \T\lo1 (g), \ldots, \T\lo k (g))
$
defines a statistical functional on $\sten G$ to $\real{}$. We call
such functionals \emph{composite linear functionals}, or C-functionals.
For example, if
\begin{eqnarray*}
T\lo1 (g) &=& \tint y g(y) \,d \me Y (y),\qquad T\lo2 (g) = \tint y\udex2 g(y) \,d \me Y
(y), \\
\rho(s\lo1, s\lo2) &=& s\lo2 - s \lo1 \udex2,
\end{eqnarray*}
%
then
$
C(g) = \var\lo g ( Y)
$
is the variance functional. The corresponding conditional statistical
functional is defined by
\begin{eqnarray*}
C\lo{x, \theta} \bigl(r(x, \cdot, \theta )\bigr) &\equal &C \of\R\lo{x, \theta }
\inv\bigl(r (x, \cdot, \theta\lo0)\bigr)
\\
&=& \rho
\bigl[ \T\lo{1,x, \theta}\bigl(r (x, \cdot, \theta)\bigr), \ldots, \T\lo {k,x,
\theta}\bigl(r (x, \cdot, \theta)\bigr) \bigr],
\end{eqnarray*}
where $\T\lo{\ell, x, \theta}$ denotes $\T\lo\ell\of\R\lo{x,
\theta}\inv$. We will use the following notation:
%
\begin{eqnarray}\label{eq:rho ell}
\rho\lo\ell(X, \theta) &\equal& \partial\rho(s) / \partial s \lo\ell
\at{s} { (T\lo{1,x, \theta }(0), \ldots, T \lo{\ell, x, \theta} (0))},
\nonumber\\
G(X, \theta ) &\equal&\bigl( \rho\lo1 ( X, \theta), \ldots, \rho\lo k (X, \theta)
\bigr)\trans,
\\
F (Y) &\equal& \bigl(f \lo1 (Y), \ldots, f \lo k (Y)\bigr)\trans. %
\nonumber\end{eqnarray}
Also note that, in our case, $ T \lo{\ell, x, \theta}(0) = \E\lo
\theta( f \lo\ell(Y) | X )\lo x$. Again, we use symbols
such as $T \lo{\ell, x}$ and $C \lo x$ to indicate $T \lo{\ell, x,
\theta\lo0}$ and $C \lo{x, \theta\lo0}$.

\begin{theorem}\label{theorem:e score l functional} Suppose the
conditions 1, 2, 3 in Theorem~\ref{theorem:general e score} hold for $C
\lo{x, \theta}$.
Then the efficient score for $\cts C$
is given by (\ref{eq:projected score}), in which
%
\begin{eqnarray}\label{eq:theorem5}
\tau\lo c (x, y, \theta\lo0) &\equal& G\trans(x, \theta\lo0) F (y ) -
G\trans(x, \theta\lo0) E \lo{\theta\lo0} \bigl( F(Y) | X \bigr) \lo x,\nonumber
\\
q \lo1 (x,\theta\lo0 ) &\equal& \partial\E\lo\theta \bigl(F\trans(Y) | X\bigr) \lo
x / \partial \theta \at\theta{ \theta\lo0} G (x, \theta\lo0 ),
\\
q \lo 2 (x,\theta\lo0) &\equal& {G\trans(x, \theta\lo0) \var \lo {\theta\lo0}
\bigl( F (Y) | X\bigr)\lo x G (x, \theta\lo0)}. \nonumber%
\end{eqnarray}
\end{theorem}

\begin{pf}
As shown in Section~\ref{section:l functional}, the Riesz
representation of $\dot\T\lo{\ell, x}(0)$ is simply~$f \lo\ell$. By
the chain rule of \frechet differentiation and definition (\ref{eq:rho ell}), we have
\[
\dot C \lo x (0) = \tsum \ell1 k \rho\lo\ell(X, \theta\lo0) \dot\T\lo{\ell, x}
(0).
\]
Hence, the Riesz representation of $\dot C \lo x (0)$ is
\[
\tau(x, y, \theta\lo0) = \tsum\ell1 k \rho\lo\ell(x, \theta\lo 0) f \lo\ell(y) =
G\trans(x, \theta\lo0) F (y).
\]
In the meantime for each $\theta$, we have
\[
\tint\tau\lo c (x, y, \theta ) \eta\lo0 (x, y, \theta ) \,d \me Y (y) = 0.
\]
Differentiate both sides of this equation with respect to $\theta$, to obtain
\[
\tint\partial\tau\lo c (x, y, \theta) / \partial\theta \eta\lo0 (x, y, \theta)
\,d \me Y (y) + \tint\tau\lo c (x, y, \theta) \putop\circ r \lo0 (x, y, \theta )
\,d \me Y (y) = 0.
\]
Hence,
\begin{eqnarray*}
&&\tint\tau\lo c (x, y, \theta) \putop\circ r \lo0 (x, y, \theta ) \,d \me Y (y)
\\
&&\qquad= - \E\lo\theta\bigl[ \partial\tau\lo c (X, Y, \theta) / \partial \theta | X
\bigr]
\\
&&\qquad = - \partial G\trans(X, \theta) / \partial\theta \bigl[ F (Y) - \E\lo\theta
\bigl(F (Y) | X\bigr) \bigr] + \partial\E\lo\theta \bigl(F\trans(Y) | X\bigr) /
\partial\theta G (X, \theta).
\end{eqnarray*}
Now take conditional expectation $E \lo\theta(\cdots| X)$ on both
sides to prove the second relation in (\ref{eq:theorem5}).
\end{pf}

It is easy to see that an alternative expression of $q \lo1 (\theta,
X)$ in Theorem~\ref{theorem:e score l functional} is
\[
q \lo1 ( x, \theta\lo0 ) = \partial\rho\bigl( E \lo\theta\bigl[f \lo1 (Y)|X
\bigr] \lo x, \ldots, E \lo\theta\bigl[ f \lo k (Y)|X\bigr] \lo x \bigr) /
\partial \theta \at\theta{\theta\lo0}.
\]
This expression is useful because $\rho( E \lo\theta[f \lo1 (Y)|X]
\lo x, \ldots, E \lo\theta[ f \lo k (Y)|X] \lo x )$
is a function of $\sigma\lo\theta(X)$, and its derivative with
respect to $\theta$ can be estimated by local linear regression,
as we will see in Section~\ref{section:estimation}.

\begin{example} For the central variance subspace, we have
\[
k = 2,\qquad f \lo1 (y) = y,\qquad f \lo2 (y) = y\udex2,\qquad
\rho(s \lo1, s \lo2) = s \lo2 - s
\lo1 \udex2.
\]
Hence, $F(y) = (y, y\hi2)\trans$, and
\begin{eqnarray*}
\rho\lo1 ( X, \theta\lo0) &\equal& \partial\bigl(s \lo2 - s \lo1 \udex 2\bigr) /\partial s
\lo1 \at{s\lo1} { E(F(Y)|X) } = - 2 \E( Y|X),
\\
\rho\lo2 (X, \theta\lo0 ) &\equal&\partial\bigl(s \lo2 - s \lo1 \udex 2\bigr) / \partial s
\lo2 \at{s\lo2} { E(F(Y)|X) } = 1.
\end{eqnarray*}
The Riesz representation of $ \dot C \lo x (0)$ and its centered
version are
\begin{eqnarray*}
\tau(X, Y, \theta\lo0)& \equal& - 2 \E( Y |X ) Y + Y\udex2,
\\
\tau\lo c (X, Y, \theta\lo0) &\equal& - 2 \E( Y |X ) Y + Y\udex2 - \E\bigl[ - 2
\E( Y |X ) Y + Y\udex2 |X \bigr]
\\
&\equal&\bigl[Y - \E(Y|X)\bigr]\udex2- \E\bigl[ \var (Y|X) | \sigma\lo{\theta\lo
0} (X)\bigr].
\end{eqnarray*}
Hence,
\begin{eqnarray*}
q \lo1 (X,\theta\lo0 )& \equal& - \partial\E\lo\theta \bigl[\tau\lo c (X, Y,
\theta)\bigr] / \partial\theta \at\theta{ \theta\lo0} = \partial\E\lo\theta \bigl[
\var\lo\theta ( Y | X) | \sigma\lo \theta (X) \bigr] / \partial\theta \at\theta {
\theta\lo0},
\\
q \lo2 (X,\theta\lo0 )& \equal& \var\bigl[ \tau\lo c ( X, Y, \theta\lo 0) | X
\bigr]= \v \bigl\{ \bigl[Y - E(Y|X)\bigr]\udex2 |X \bigr\}.
\end{eqnarray*}
In this case, the subspace $\sten T \lo\phi\oc$ consists of
functions of the form
\[
\bigl[ h(x) - E \bigl(h(X) |\sigma\lo{\theta\lo0} (X)\bigr)\bigr] \bigl[\bigl(Y
- \E(Y|X)\bigr)\udex 2- \E\bigl( \var (Y|X) | \sigma\lo{\theta\lo0} (X)\bigr)
\bigr],
\]
where $h$ is an arbitrary member of $\lt{\lambda\lo0}$.
Interestingly, the estimating equation proposed by Zhu, Dong and Li
(\citeyear{ZhuDonLi13})
is a special member of $\sten T \lo\phi\oc$ with $h(x)$ taken to be
the components of $x$. 
\end{example}

\section{Implicit statistical functionals}\label{section:i functional}

\def\ka{\kappa}
\def\I{I}

In this section, we study statistical functionals defined implicitly
through estimating equations. Many robust estimators, such
as conditional medians and quantiles, are of this type. Let $\Xi
\subseteq\real{}$,
and $e\dvtx \Xi\times\om\lo Y \to\real{}$ be a function of the
parameter $\xi$ and the variable $y$. Such functions
are called estimating functions [see, e.g., \citet{God60}, \citet{LiMcC94}]. If
the equation
%
\begin{equation}
\tint\lo{\om\lo Y} e ( \xi, y) g(y) \,d \me Y( y) = 0 \label{eq:implicit}
\end{equation}
has a unique solution for each $g \in\sten G$, then it defines a functional
$I\dvtx \sten G \to\real{}$ that assigns each $g$ the solution to (\ref
{eq:implicit}).
We call such functionals \emph{implicit functionals}, or I-functionals.
If we replace $g \in\sten G$ by a conditional density function $\eta
(\cdot, \cdot, \theta) \in\sten H$, then (\ref{eq:implicit}) becomes
\[
\tint\lo{\om\lo Y} e ( \xi, y) \eta(x, y, \theta) \,d \me Y( y) = 0.
\]
The corresponding conditional statistical functional is $\I\lo{x,
\theta} ( r (x, \cdot, \theta)) =\break   \I\of\R\lo{x, \theta}\inv(r(x,
\cdot, \theta))$.
The $I$-central subspace is defined by the statement
\[
\I\lo X (0) \mbox{ is measurable with respect to } \sigma\lo{\theta
\lo0} (X).
\]
Naturally, we write the function $(\theta, x) \mapsto I \lo{x, \theta}
(0)$ as $\xi( \theta, x )$.

To simplify the presentation, we use the notion of generalized
functions. Let $\sten K$
be the class of functions defined on a bounded set $B$ in ${\mathbb R}$
that have derivatives of all orders, whose topology is defined
by the uniform convergence of all derivatives. Any continuous linear functional
$U\dvtx \sten K \to\real{}$ with respect to this topology is called a
generalized function. For example, let $a \in B$. Then
it can be shown that the linear functional
\[
\delta\lo a\dvtx \sten K \to\real{},\qquad \phi\mapsto\phi(a)
\]
is continuous with respect to this topology. This continuous linear
functional is called the Dirac delta function. We identify the
functional $\delta\lo a$
with an imagined function $x \mapsto\delta\lo a (x)$ on $B$ and write
$\delta\lo a(\phi)$
formally as the integral
\[
\delta\lo a (\phi) = \tint\phi(x) \delta\lo a (x) \,d \lambda(x).
\]
A consequence of this convention is that if we
pretend $\delta\lo a (x) $ to be the derivative $\partial I ( x \le a)
/ \partial a$ of the indicator function $I (x \le a)$ then
we get correct answers at the integral level. For example, for any
constant $a$ and small number $\epsilon$, we have
\[
\tint\bigl[ I ( y \le a + \epsilon) - I (y \le a) \bigr] g(y) \,d y = \int\delta
\lo a (y) \epsilon g (y) \,d y + o (\epsilon) = \epsilon g(a) + o(\epsilon).
\]
Thus, the pretended linearization $I ( y \le a + \epsilon) - I (y \le
a) = \delta\lo a(y) \epsilon+ o(\epsilon)$ has caused no inconsistency.
We use this device to simplify
our presentation of quantiles.

\begin{theorem}\label{theorem:i-functional} Suppose the conditions 1,
2, 3 in Theorem~\ref{theorem:general e score} hold for $I \lo{x, \theta
}$. Moreover,
suppose
for each $\xi\in\Xi$, $g \in\sten G$, there is a (generalized)
function $\dot e (\xi, y)$,
which plays the role of $\partial e ( \xi, y) / \partial\xi$, such that
%
\begin{equation}
\biggl\llvert \tint\lo{\om\lo Y} \bigl[ e (\xi+ a, y) - e (\xi, y) - \dot e (\xi
, y) a \bigr] g(y) \,d \me Y (y) \biggr\rrvert = o(a). \label{eq:generalized function}
\end{equation}
Then the efficient score for the $I$-central subspace is (\ref
{eq:projected score}) in which
%
\begin{eqnarray}
\label{eq:q1 q2 for i functional} %
\tau\lo c (x, y, \theta\lo0 )&\equal& - { e \bigl(\xi(
\theta\lo0, x), y\bigr)}/E \lo{\theta\lo0}\bigl[\dot e \bigl(\xi(\theta\lo0, x),
Y\bigr)|X \bigr]\lo x,
\nonumber\\
q \lo1 (x,\theta\lo0) &\equal& \partial\xi(\theta,x) / \partial \theta \at\theta{
\theta\lo0},
\\
q \lo2 (x,\theta\lo0 )& \equal& {\E\lo{\theta\lo0}\bigl[ e\udex2 \bigl( \xi(
\theta, X),Y\bigr)|X \bigr]\lo x}/\E\lo{\theta\lo0}\udex2 \bigl[\dot e \bigl(\xi
(\theta\lo0, X), Y\bigr) |X \bigr] \lo x. \nonumber%
\end{eqnarray}
\end{theorem}

\begin{pf} Differentiating
both sides of the equation
\[
\tint\lo{\om\lo Y} e \bigl( I \lo{x, \theta} \bigl(\epsilon r(x, \cdot, \theta )
\bigr), y\bigr) \bigl( 1 + \epsilon r(x, y, \theta) / 2\bigr)\udex2 \eta\lo0 (x,
y, \theta) \,d \me Y( y) =0
\]
with respect to $\epsilon$ at $\epsilon= 0$, and using the relation
\[
\partial\bigl( 1 + \epsilon r(x, y, \theta) / 2\bigr)\udex2 / \partial \epsilon
\at\epsilon0 = 2 \bigl( 1 + 0 r(x, y, \theta)/2\bigr) r(x,y,\theta)/2 = r(x,y,
\theta),
\]
we find
\[
\partial\I\lo{x, \theta} \bigl(\epsilon r(x, \cdot, \theta) \bigr) / \partial
\epsilon\at\epsilon0 = - \E \biggl[ \frac{ e ( \I\lo{x, \theta} (0), Y) }{\E(\dot e ( \I
\lo{x, \theta} (0), Y) |X )\lo x} r(x,Y,\theta) |X \biggr]\lo
x.
\]
Hence, the Riesz representation of the \frechet derivative $\dot\I\lo
{x, \theta} (0)$ is
%
\begin{eqnarray}
\label{eq:tau and tau c} %
\tau(x,y,\theta)& \equal& - { e \bigl( \I\lo{x, \theta}
(0), y\bigr)}/\E\bigl[\dot e \bigl( \I\lo{x, \theta} (0), Y\bigr) |X \bigr]\lo x
\nonumber
\\[-8pt]
\\[-8pt]
\nonumber
&\equal& - { e \bigl(\xi(\theta, x), y\bigr)}/E\bigl[\dot e \bigl(\xi(\theta, x)
, Y\bigr)|X \bigr]\lo x=\tau\lo c (x, y, \theta), %
\end{eqnarray}
where the last equality holds because, by definition, $\E[ e (\xi
(\theta, X), Y ) | X ] = 0$.
By~(\ref{eq:tau and tau c}) and the definition of $q \lo1$ in
Theorem~\ref{theorem:general e score}, 
\begin{equation}\qquad
q \lo1 (x, \theta\lo0) = - \E\bigl[ e \bigl(\xi (\theta\lo0, x), y\bigr)
\score x y | X \bigr] \lo x/\E\bigl[\dot e \bigl(\xi(\theta\lo0, x), Y\bigr)|X
\bigr]\lo x. \label{eq:this is equal to}
\end{equation}
To further simplify the numerator of the right-hand side, differentiate
both sides of the equation
$
\tint e (\xi(\theta,x), y) \eta\lo0 (x, y, \theta) \,d \me Y (y) = 0
$
to obtain
\begin{eqnarray*}
&& \biggl[ \tint\dot e \bigl(\xi(\theta\lo0, x), y\bigr) \eta\lo0 (x, y, \theta
\lo0 ) \,d \me Y (y) \biggr] \dot\xi(\theta\lo0,x)
\\
&&\qquad{} + \tint e \bigl(\xi(\theta\lo0,x), y\bigr) \score x y \eta\lo0
(x, y, \theta \lo0 ) \,d \me Y (y) = 0,
\end{eqnarray*}
where $\dot\xi(\theta\lo0, x)$ denotes $\partial\xi(\theta, x)
/ \partial\theta \at\theta{ \theta\lo0}$.
Hence,
\[
\E\bigl[ e \bigl(\xi(\theta\lo0, X), Y\bigr) \score XY | X \bigr]\lo x = - \E
\bigl[\dot e \bigl(\xi(\theta\lo0,X), Y\bigr) |X \bigr]\lo x \dot\xi(\theta
\lo0,x).
\]
Substitute this into (\ref{eq:this is equal to}) to prove the first
relation in (\ref{eq:q1 q2 for i functional}).
Substitute (\ref{eq:tau and tau c}) into the definition of $q \lo2$ in
Theorem~\ref{theorem:general e score} to prove the second relation
in (\ref{eq:q1 q2 for i functional}).
\end{pf}

In the next example, we derive the efficient score for a particular
type of I-functional---the quantile functional.

\begin{example} If we assume all densities in $\sten G$ are continuous,
then the $p$th quantile is the solution to the equation
$
\E\I( Y \le\xi) = p
$. Equivalently, $\xi$ is the solution to
the equation
\[
E\bigl[ e (\xi, y ) \bigr] = E\bigl[ - \sgn(y - \xi) + 1 - 2p \bigr] = 0.
\]
Because the generalized derivative of $\sgn(t)$ is $2 \delta\lo0
(t)$, we have
$
\dot e ( \xi, y) = 2 \delta\lo\xi(y)$.
Hence
\[
\E \bigl[\dot e \bigl(\xi(\theta, X), Y\bigr) |X \bigr]\lo x \equal \tint\lo{
\om\lo Y} 2 \delta\lo{\xi( \theta, x)} (y ) \eta\lo0 (x, y, \theta) \,d y = 2 \eta
\lo0 \bigl(x, \xi(\theta,x), \theta\bigr).
\]
Because $e ( \xi(\theta, x), Y)$ is a binary random variable that
takes the value $2 (1-p)$ with probability $p$
and $-2p$ with probability $(1-p)$, we have
\[
\E\bigl[ e\udex2 \bigl( \xi(\theta, x), Y \bigr) | X \bigr]\lo x \equal
\bigl[2(1-p)\bigr]\udex2 p + (-2p)\udex2 (1-p) = 4 (1-p) p.
\]
Hence, in the efficient score,
\[
q \lo2 (x,\theta\lo0 ) = \eta\udex2 \lo0 \bigl(x, \xi(\theta\lo0, x), \theta\lo0
\bigr)/ \bigl[(1-p)p\bigr],
\]
and $q \lo1 (x, \theta\lo0)$ is as given in (\ref{eq:q1 q2 for i
functional}) with $\xi(\theta, x)$ being the conditional $p$th
quantile of~$Y$ given $X$.
\end{example}

\section{Effect of estimating the central subspace}\label{section:effect}

Throughout the previous sections, we have treated $X$ as the true
predictor from the central subspace; that is,
$X = \zeta \trans \tilde X$ where $\tilde X$ is the original
predictor and $\zeta$ is a basis matrix of the central subspace
based for $Y | \tilde X$. However, in practice, $\zeta$ itself needs to
be estimated and, in theory at least, should affect
the form of the efficient score about $\beta$. While our simulation
studies in Section~\ref{section:simulation} indicate that this effect
is very small,\vadjust{\goodbreak}
for theoretical rigor we present here the efficient score treating the
central subspace as an additional (finite-dimensional) nuisance parameter.
For convenience, we use $\zeta$ to denote both a $p \times d$ matrix
and the corresponding $(p-d)d$-dimensional Grassmann manifold.

Let $S\eff(\zeta\trans\tilde X, Y, \theta)$ denote the efficient
score in Theorem~\ref{theorem:general e score} with $X$ replaced by
$\zeta\trans\tilde X$.
Let
$S\eff\udex* (\tilde X, Y, \zeta, \theta)$ denote the efficient score
for $\theta$ with $\zeta$ treated as an additional
nuisance parameter.
Let $s\lo\zeta(\tilde x, y, \zeta, \theta)$ denote the score with
respect to~$\zeta$. In the supplementary material [\citet{supp}] (Section \abb{iii})
it is shown that
%
\begin{equation}
S\eff\udex* = S\eff - \Pi\bigl(S\eff | \operatorname{span} \bigl(\Pi\bigl(s\lo \zeta|
\sten T\lo\phi\oc\bigr) \bigr) \bigr) \equiv S \eff- g, \label{lemma:goal}
\end{equation}
where $g$ is the function
\begin{eqnarray*}
g &=& q\lo2\inv\bigl(\zeta\lo0 \trans\tilde x, \theta\lo0\bigr) \tau\lo c \bigl(\zeta\lo0\trans
\tilde x, y, \theta\lo0\bigr) E \bigl\{q\lo3 \bigl(\zeta\lo0 \trans\tilde X, \theta\lo0\bigr) q
\lo3\udex{*{\mathsf T}}\bigl(\zeta\lo0 \trans\tilde X, \theta\lo0\bigr) q\lo 2\inv\bigl(
\zeta\lo0 \trans\tilde x,\theta\lo0\bigr)\bigr\}
\\
&&{}\times E \bigl\{qe\udex* \lo3\bigl(\zeta\lo0\trans\tilde X, \theta\lo0\bigr) q\udex {*
{\mathsf T}} \lo3 \bigl(\zeta\lo0\trans\tilde X, \theta\lo0\bigr) q\lo2\inv\bigl(\zeta\lo0
\trans\tilde x,\theta\lo0\bigr)\bigr\}\udex\dagger q\udex* \lo3 \bigl(\zeta\lo0\trans\tilde
x, \theta\lo0\bigr),
\end{eqnarray*}
in which $\dagger$ indicates the Moore--Penrose inverse, and
\begin{eqnarray*}
q \udex*\lo1 \bigl(\zeta\lo0\trans\tilde x, \theta\lo0\bigr) &=& \E \bigl[s\lo \zeta(
\tilde X, Y, \zeta\lo0, \theta\lo0 ) \tau\lo c \bigl(\zeta\lo0 \trans\tilde X, Y,
\theta\lo0\bigr) | \zeta\lo0\trans\tilde X\bigr] \lo {\tilde x}
\\
q \udex* \lo3\bigl (\zeta\lo0\trans\tilde x, \theta\lo0\bigr) &=& q \udex* \lo1\bigl (\zeta\lo0
\trans\tilde x, \theta\lo0\bigr)
\\
&&{} - \E\lo{\theta\lo0} \bigl[q \udex*\lo1 \bigl(\zeta\lo0 \trans\tilde X, \theta\lo0\bigr)
q \lo2 \inv \bigl(\zeta\lo0 \trans\tilde x, \theta\lo0\bigr) | \sigma\lo{\theta\lo0} \bigl(\zeta
\lo0\trans\tilde X\bigr)\bigr] \lo {\tilde x}\\
&&\quad{}/\E\bigl[q \lo2 \inv \bigl(\zeta\lo0 \trans
\tilde x, \theta\lo 0\bigr) | \sigma\lo{\theta\lo0} \bigl(\zeta\lo0 \trans\tilde X\bigr)
\bigr]_{\nano
\tilde x}.
\end{eqnarray*}
In theory, the asymptotic variance bound based on $S\eff$ is lower than
or equal to that based on $S\eff\hi*$. However, in the simulation studies
(Table~\ref{tab2} in Section~\ref{section:simulation}) we see that the
theoretical lower bound based on $S\eff$ is nearly reached, which
indicates that effect of estimating
the central subspace on the efficient score for $\beta$ is very small.

\section{Estimation}\label{section:estimation}

In this section, we introduce semiparametrically efficient estimators
using the theory developed
in the previous sections. For the L-functionals, we develop the
estimator in full generality, but
for the C- and I-functionals we focus on the conditional variance
functional and the conditional quantile functional.
Procedures for other C- or I-functionals can be developed by analogy.

\def\dmave{{ {DMAVE}}}

\def\aaa{{ {A}}A{{A}}}

We first clarify two points related to the algorithm we will propose.
First, since we will rely heavily on
the \mave-type algorithms, it is more convenient to use the $\beta
$-parameterization rather than the $\theta$-parameterization,
and avoid redundancy in $\beta$
by taking the generalized matrix inverse. We will justify the $\beta
$-parameterization after introducing the algorithm.
Second, the \mave\ algorithm actually has two variants: the outer
product gradient \mbox{(\opg)} and a refined version of \mave\ (\rmave). Typically,
\opg,  \mave\ and \rmave\ are progressively more accurate and
require more computation. In the following, the \mave-type algorithm
can be replaced either by \rmave\ for greater accuracy or by \opg\
for less computation.

\def\see{{{SEE}}}

Our estimation procedure is divided into four steps. In step 1, we
estimate the central subspace and project $\tilde X$ on to this
subspace to obtain $X$.
In step 2, we estimate
$T(\eta\lo0 ( X \lo1, \cdot)), \ldots, T(\eta\lo0( X \lo n, \cdot
))$ using a $d$-dimensional kernel
estimate. These estimates are used as the proxy response, and we denote
them by $\hat Y \lo1, \ldots, \hat Y \lo n$. In
step 3, we apply \mave\  to $( X \lo1, \hat Y \lo1 ), \ldots, ( X
\lo n, \hat Y \lo n)$ to estimate
an initial value for~$\beta$. In
step 4, we use one-step Newton--Raphson algorithm based on the
efficient score and efficient information
to approximate the semiparametrically efficient estimate.
We call our estimator \see, which stands for \emph{semiparametrically
efficient estimator}.

\textit{Preparation step: A \mave\ code}. Let $(X \lo1, U \lo1),
\ldots, (X \lo n, U \lo n)$ be a random sample from $(X, U)$, and
$K \lo h(\cdot)$ be a kernel with bandwith $h$. That is, $K \lo h (t)
= K (t / h)/h$ for some symmetric function $K$ that integrates to 1 and
$h > 0$.
Compute the objective function
%
\begin{equation}
\label{eq:general objective} \Gamma(a, b, A) = \tsum i 1 n \tsum j 1 n
\bigl[ U \lo j - a \lo
i - b \lo i\trans A\trans(X \lo j - X \lo j ) \bigr] \hi2 K \lo h ( X \lo j - X
\lo i),
\end{equation}
where
$a \lo1, \ldots, a \lo n \in\real{}$, $b \lo1, \ldots, b \lo n \in
\real {d}$, $a$ and $b$ on the left denote
$ (a \lo1, \ldots, a \lo n)\trans$ and $ (b \lo1\trans, \ldots, b
\lo n\trans)\trans$, respectively, and $A \in\real{p \times d}$. We
will use
this objective function in several ways.
It is well known that minimization of (\ref{eq:general objective}) can
be solved
by iterations of least squares, and in each iteration there is an
explicit solution.
Thus, purely search-based numerical optimization (such as the simplex
method) is avoided.
See, for example, \citet{LiLiZhu10} and \citet{YinLi11}.

\textit{Step \textup{1:} Estimation of central subspace.}
We use the \mave-ensemble in \citet{YinLi11}
to estimate the central subspace. In this procedure, $U$ in (\ref
{eq:general objective}) is taken to be a set of functions
$\{f\lo1(Y),\ldots,f\lo m(Y)\}$ randomly sampled from a dense family
in $L \lo2 ( \mu\lo Y)$. In this paper, take this set to be $\{(\sin
(t \lo i y), \cos(t \lo i y))\dvtx i = 1,\ldots, 10\}$, where $t \lo1,
\ldots, t \lo{10}$
are i.i.d. $\operatorname{unif}(0,4)$. The sample of responses $Y\lo1, \ldots,
Y \lo n$ are standardized so that
the range $(0,4)$ of the uniform distribution represents a reasonably
rich class of functions relative to the range of $Y$. The basis matrix
$\zeta$ of $\ca S \lo{Y|X}$ is
then estimated by the \mave-ensemble. The projected predictor $X = \hat
\zeta\trans\tilde X $ is taken as the predictor in steps 2--4.
Since our goal is to estimate $\cts T$, the choice of the working
dimension $\hat d$ of $\cs$ is not crucial.

As was
shown in \citet{YinLi11}, at the population level, the \abb{mave}
ensemble is guaranteed to recover central subspace exhaustively
as long as the functions of $Y$ form a characterizing family.
In practice, it is true any information lost in the initial step will
be inherited by \see. However, our experiences indicate that
this problem can be mitigated by
using a sufficiently rich ensemble family---for example, by
increasing the range of the uniform distribution and the number of $t
\lo i$'s.

Several other methods are available for exhaustive estimation of the
central subspace, such as the
semiparametrically efficient estimator of \citet{MaZhu13N1},
the \abb{dmave}
by \citet{Xia07} and the Sliced Regression by \citet
{WanXia08}. Here, we have chosen the \abb{mave}-ensemble for its
computational simplicity.

\textit{Step \textup{2:} Estimation of proxy response.}
This step is unnecessary for the L-functionals: because $T ( \eta( X,
\cdot) ) = E [ f(Y) | X ]$ for some function $f$,
we can use $f(Y)$ itself as the proxy response $\hat Y$.
For the conditional variance functionals, this step needs not be fully
performed: we can use $ (Y - \hat E (Y| X ))\hi2$ as
the proxy response $\hat Y$, where $\hat E (Y|X )$ is the kernel estimator
of $E(Y|X )$. If simplification of this type is not applicable, then we
need to perform nonparametric estimation of $T(\eta(X, \cdot))$. For example,
for the I-functionals, we use the minimizer $\xi\lo i \udex*$ of the function
$
E \lo n [e ( Y, \xi) K \lo h ( X - X \lo i)]
$
as the proxy
response $\hat Y \lo i$.

\textit{Step \textup{3:} Initial estimate of $\beta$}. Apply \mave\ to $(X
\lo1, \hat Y \lo1), \ldots, (X \lo n, \hat Y \lo n)$
to obtain an initial estimate of $\beta$ by minimizing $\Gamma(a, b,
\beta)$ over all $a, b, \beta$. Denote this initial estimate as $\tilde
\beta$.

\def\eff{_{\mathrm{ eff}}}

\textit{Step \textup{4:} The one-step Newton--Raphson algorithm.}
Rather than attempting to solve the score equation (\ref{eq:score
equation}), we propose a one-step Newton--Raphson
procedure. Let $S\eff(X, Y, \beta)$ be the efficient score in $\beta$,
which is obtained by
replacing
$\sigma\lo\theta(X)$ by $\beta\trans X$ and $\partial/ \partial
\theta$ by $\partial/ \partial\vec(\beta)$
whenever applicable. Let
$J \neff(\beta)=E \lo n [ S \eff(X, Y, \beta) S_{\mathrm{\scriptsize
eff}} \trans (X, Y, \beta)]$.
We estimate $\beta$ by
%
\begin{equation}
\label{eq:Newton-est beta} \hat\beta= \tilde\beta+ \bigl[J\neff(\tilde\beta)\bigr]\udex
\dagger E \lo n \bigl[ S\eff( X, Y, \tilde\beta) \bigr],
\end{equation}
where
$\tilde\beta$ is the initial value from step 3.

\def\mave{{{{MAVE}}}}

We now describe in detail how to compute $E \lo n [ S\eff( X, Y, \tilde
\beta) ]$ for different functionals. For the L-functional,
the efficient score involves the following functions:
\[
\textup{(a)}\quad E\bigl[ f(Y) | \beta\trans X\bigr],\qquad \textup{(b)}\quad E \bigl[ f\udex2(Y) | X
\bigr],\qquad \textup{(c)}\quad \partial E \bigl[ f(Y) | \beta\trans X\bigr]/ \partial\vec(\beta).
\]
For the C-functional, it involves the functions:
\begin{eqnarray*}
&&\mathrm{(d)}\quad E \lo{\theta\lo0} \bigl[ f \lo\ell(Y) | X \bigr],\qquad
 \mathrm{(e)}\quad \partial
\rho\bigl( E \lo\beta\bigl[ f \lo1 (Y) |X \bigr], \ldots, E \lo\beta\bigl[ f \lo
\ell(Y) | X \bigr] \bigr) / \partial\vec(\beta),
\\
& &\mathrm{(f)}\quad \cov \bigl[f \lo i (Y), f \lo j (Y) | X \bigr]. %
\end{eqnarray*}
For the conditional quantile functional, it involves the functions:
\begin{eqnarray*}
&&\mathrm{(g)}\quad E\bigl[ \xi( X ) | \beta\trans X \bigr],
\qquad \mathrm{(h)}\quad \partial E\bigl[ \xi( X ) |
\beta \trans X \bigr] / \partial\vec(\beta),\\
&&\mathrm{(i)} \quad\eta\lo0 \bigl( X, E\bigl[ \xi( X
) | \beta\trans X \bigr]\bigr).
\end{eqnarray*}
These functions can be categorized into three types: (a), (b), (d),
(f), (g) 
are conditional means of random variables;
(c), (e), (h)
are derivatives of functions of $\beta\trans X$ with respect to $\vec
(\beta)$;
(i) is the conditional density evaluated at a quantile.
The first two types can be solved by minimizing $\Gamma(a, b, A)$ in
(\ref{eq:general objective}) with specific $U\lo i$ and~$A$:
Table \ref{tab1} gives the details of these random variables and matrices, as well
as which parts of the \mave\ output are needed in \see.

\begin{table}
\caption{Using \mave\ to estimate quantities \textup{(a)} through \textup{(h)}
in efficient score}\label{tab1}
\begin{tabular*}{\textwidth}{@{\extracolsep{\fill}}lccc@{}}
\hline
\textbf{Quantities} & $\bolds{A}$ & $\bolds{U\lo i}$ & \multicolumn{1}{c@{}}{\textbf{MAVE output}} \\
\hline
(a) & $\beta$ & $f(Y)$ & $a\udex*$ \\
(b) & $I \lo p$ & $f\udex2(Y)$ & $a\udex*$ \\
(c) & $\beta$ & $f(Y)$ & $b \udex*$ \\
(d) & $I \lo p$ & $f \lo\ell(Y)$ & $a \udex*$ \\
(e) & $\beta$ & $\rho( E \lo n [ f \lo1 (Y)|X], \ldots, E \lo n [ f
\lo k (Y) | X ] )$ & $b \udex*$ \\
(f) & $I \lo p$ & $f \lo\ell(Y)$, $f \lo\ell(Y) f \lo{\ell'} (Y)$
& $a \udex*$ \\
(g) & $\beta$ & $\xi\lo n ( X)$ & $a \udex*$ \\
(h) & $\beta$ & $\xi\lo n ( X)$ & $b \udex*$ \\
\hline
\end{tabular*}
\end{table}

Finally, we estimate $\eta\lo0 (x, y \lo0)$ for any $y \lo0$ by the
kernel conditional density estimator:
\[
E \lo n \bigl[ K \lo{h \lo1} ( Y - y \lo0) K \lo{h} ( X - x) \bigr] / E \lo n
\bigl[ K \lo{h} ( X - x ) \bigr],
\]
where $h \lo1$ is a different bandwidth for the response variable.

\textit{Tuning of bandwidths}. We need to determine the kernel
bandwidths $h$ at various stages in the above algorithm.
We use the Gaussian kernel with optimal bandwidth $h =c n \udex
{-1/(p+4)}$ for nonparametric regression
[see, e.g., \citet{Xiaetal02}], where $c$ is determined by five-fold
cross validation.
That is, we randomly divide the data into 5 subsets of roughly equal
sizes. For each $i = 1, \ldots, 5$, we use the $i$th subset as the
testing set and the rest as the training set. For a given $c$, we
conduct dimension reduction on the training set,
and use the sufficient predictor to evaluate a certain prediction
criterion at each point on the testing set and average these
evaluations over the testing set, and finally average the five averages
of the criterion to obtain a single number. We then minimize the resulting
criterion over $c$ by a grid search. Naturally, the prediction
criterion depends on the object to be estimated using that kernel.
Below is a list of the prediction criteria we propose for the four
steps in the estimation procedure.

In step 1: We use the distance correlation introduced by Sz\'{e}kely and
Rizzo [(\citeyear{SzeRiz09}), Theorem~1, expression
(2.11); $p$ and $q$ therein
are taken to be 2].

In step 2: For the L-functionals, no tuning is needed. For the
conditional variance functional, we need to estimate $E(Y|X)$, and
we use the prediction criterion
$[Y-\hat E (Y|X)]\udex2$, where $\hat E (Y|X)$ is the kernel estimate
of $E(Y|X)$ based on the training set
using a tuning constant $c$. For conditional median,
we use the prediction criterion $| Y - \hat M (Y|X) |$, where $\hat M
(Y|X)$ is the kernel estimate of the conditional median based on the
training set using a specific tuning constant.

In step 3: We use the prediction criterion $[\hat Y - \hat E ( \hat Y |
\tilde\beta\trans X ) ]\udex2$, where $\hat Y$ is
the proxy response obtained from step 2, and $\hat E ( \hat Y | \tilde
\beta\trans X ) $ is the \mave\ output $a\udex*$.

In step 4: There are three types of kernels in this step: the kernel
for $X$, the kernel for $\tilde\beta\trans X$, and the kernel for $Y$
(the last one is needed only for the conditional median functional).
Corresponding to these types, we use bandwidths
\[
h = c n\udex{ - 1/(d + 4)}, \qquad h = c n \udex{ - 1 / (s + 4)},\qquad  h = c n \udex{ - 1/
(1+4)}.
\]
We use cross validation to determine the common $c$. Once again,
we use different prediction criteria for different functionals. For the
conditional mean functional, we use
the criterion $[Y - \hat E(Y | \hat\beta\trans X ) ]\udex2$.
For conditional variance
functional, we use
the criterion $\{(Y - \hat E (Y|X) )\udex2 - \hat E [ Y - \hat E (Y|X)
| \hat\beta\trans X ]\udex2 \}\udex2$.
For the conditional median functional, we use the criterion\break
$| Y - \hat E [ \hat M (Y|X) | \hat\beta\trans X ] |$.

\def\mpinv{^{\nano\dagger}}

\textit{Justification of parameterization.}
We now justify the one-step iteration formula~(\ref{eq:Newton-est
beta}) as an equivalent form of
%
\begin{equation}
\label{eq:one step in theta} \hat\theta= \tilde\theta+ J\neff(\tilde\theta)\inv E \lo n \bigl[S
\eff (X, Y, \tilde\theta)\bigr].
\end{equation}
Since $\beta$ is a function of a $s(d-s)$ dimensional parameter $\theta
$, the efficient information $J\neff(\beta)$ has rank
$s(d-s)$. Let $\Gamma$ denote the $sd \times[s(d-s)]$ matrix whose
columns are eigenvectors of $J\neff(\beta)$ corresponding to its nonzero
eigenvalues, and let $\beta= \Gamma\theta$, where $\theta$ is a free
parameter in $\real{s(d-s)}$. In this parameterization,
\[
S \eff(X, Y, \theta) = \Gamma\trans S \eff(X, Y, \beta),\qquad J \neff( \theta) =
\Gamma\trans J \neff(\beta) \Gamma.
\]
Hence, (\ref{eq:one step in theta}) is equivalent to
\[
\hat\theta= \tilde\theta+ \bigl[\Gamma\trans J \neff(\tilde\beta) \Gamma\bigr]
\inv\Gamma\trans E \lo n \bigl[S\eff(X, Y, \tilde\beta)\bigr].
\]
Multiply both sides by $\Gamma$ from the left, we find
\[
\Gamma\hat\theta= \Gamma\tilde\theta+ \Gamma\bigl[\Gamma\trans J \neff(\tilde
\beta) \Gamma\bigr]\inv\Gamma\trans E \lo n \bigl[S\eff(X, Y, \tilde\beta)\bigr].
\]
By construction, $\Gamma\hat\theta= \hat\beta$, $\Gamma\tilde
\theta= \tilde\beta$, and $\Gamma[\Gamma\trans J \neff(\tilde\beta
) \Gamma]\inv\Gamma\trans$
is Moore--Penrose inverse of $J\neff(\tilde\beta)$. Thus, the above
iterative formula is the same as (\ref{eq:Newton-est beta}).

In concluding this section, we point out two attractive features of our
algorithm, which were briefly touched on in the \hyperref[sec1]{Introduction}.
First, since our algorithm is implemented by repeated applications of
variations of \mave, it essentially consists of sequence of
least squares algorithms, thus avoiding any search-based numerical
optimization, which can be infeasible when the dimension of $\theta$ is
high. The second advantage is that, since our algorithm is based on the
$\beta$-parameterization, we do not need any subjectively chosen
parameterization. In comparison, \citet{MaZhu13N1} used the
parameterization $\beta= (I \lo s, \theta)\trans$, where $\theta\in
\real{(d - s) \times s}$
is a matrix with free-varying entries. Note that this is not without
loss of generality, because in reality the first
$d$ rows of $\beta$ can be linearly dependent.

\section{Simulation comparisons}\label{section:simulation}

In this section, we conduct simulation comparisons between \see\  and
other methods for estimating three types of $T$-central subspaces:
conditional mean, conditional variance and conditional quantile. We use
the distance between two subspaces proposed by \citet{LiZhaChi05}
to measure estimation errors, which is defined as
%
\begin{equation}
\mbox{dist}(\ca S \lo1, \ca S \lo2) = \|\Pi\lo{\ca S \lo1} - \Pi\lo {\ca S \lo2}
\| \lo2, \label{eq:lzc}
\end{equation}
where $\ca S \lo1$ and $\ca S \lo2$ are subspaces of $\real p$, and
$\|\cdot\|\lo2$ is the $L\lo2$ norm in $\real{p \times p}$. For each
of the following models,
the sample size is taken to be $ n = 200$ or 500, or both; each sample
is repeatedly drawn for $n \lo{\mathrm{sim}} = 100$ times in the simulation.
In all simulations, we fix the working dimension in step 1
at $\hat d = 3$, even though $d$ is no greater than 2 in all examples.

The explicit forms of the efficient scores efficient information for
Models \abb{i}--\abb{vii} used in the following comparisons
are derived in the supplementary material [\citet{supp}] (Section \abb{ii}).

\def\model{\mbox{{{M}}odel}}
\def\fir{\mbox{{I}}}
\def\sec{\mbox{{{II}}}}
\def\thi{\mbox{{{III}}}}
\def\fou{\mbox{{{IV}}}}
\def\fiv{\mbox{{{V}}}}
\def\six{\mbox{{{VI}}}}
\def\sev{\mbox{{{VII}}}}
\def\eig{\mbox{{{VIII}}}}

(a) \textit{Comparison for the central mean subspace.}
In this case, the functional $T(\eta(X, \cdot))$ is the conditional
mean $E(Y|X)$.
We compare \see\ with \abb{rmave} under the following models:
\begin{eqnarray*}
&& \model\ \fir\dvtx\quad Y = X\lo1 + \bigl( 1 + | X\lo2 | \bigr) \varepsilon,
\\
&& \model\ \sec\dvtx\quad Y =X \lo1 ( X \lo1 + X \lo2 + 1 ) + 0.5 \varepsilon, \\
&& \model\ \thi\dvtx\quad Y = X\lo1 +\bigl( 1 + | X\lo1 | \bigr) \varepsilon,
\\
&& \model\ \fou\dvtx\quad Y|X \sim \operatorname{Poisson}\bigl(|X\lo1+X\lo2|\bigr),
\end{eqnarray*}
where $X \sim N(0, I\lo{10})$, $X \indep\varepsilon$, and $\varepsilon
\sim N(0,1)$. These models represent a variety of scenarios one might
encounter in practice.
Specifically, the central
mean subspace is a proper subspace of the central subspace in Model
$\fir$, but coincides with the latter in the other models.
The conditional variance $\v(Y | X)$ is a
constant in Model $\sec$, but depends on $X$ in the other models.
Because of its additive error structure
Model $\sec$ is favorable to \mave. Finally, Model $\fou$ has a
discrete response and the error only enters implicitly.
Model $\fir$ and
Model $\thi$ will be used again for Comparison (b), where conditional
variance is the target;
Model \sec\ was also used in \citet{Li91N2} and \citet{Xiaetal02}.

\def\nsim{n \lo{\mathrm{sim}}}

The results with sample sizes $n = 200$ and $n=500$ are presented in
Table~\ref{tab2}, in the blocks indicated by $E(Y|X)$.
The entries are in the form $a(b)$, where $a$ is the mean, and $b$ the
standard error, of the distance
(\ref{eq:lzc}) between the true and estimated $\cts E$, based on $\nsim
= 100$ simulated samples.

(b) \textit{Comparison for the central variance subspace}. Let
$T(\eta(X, \cdot))$ be the conditional variance $\var(Y|X)$.
We
compare \see\  with the estimator proposed in \citet{ZhuZhu09}
and Zhu, Dong, and Li (\citeyear{ZhuDonLi13}).
In Model $\fir$, the central variance subspace is different from either
the central mean subspace or the central subspace,
while in Model $\thi$, the three spaces coincide. The results with $n =
200$ and $n = 500$ are reported in Table~\ref{tab2}, in the blocks indicated by
$\var(Y|X)$.

%
\begin{table}
\caption{Comparison of \see\ with other estimators for three
statistical functionals}\label{tab2}
\begin{tabular*}{\textwidth}{@{\extracolsep{\fill}}lcccccc@{}}
\hline
$\bolds{n}$& {\textbf{Functionals}} & {\textbf{Models}} &
\multicolumn{4}{c@{}}{\textbf{Estimators}} \\
\hline
{200}& {$E(Y|X)$} & &\multicolumn
{2}{c}{ \rmave} & \see & ALB \\
&& $\fir$ &\multicolumn{2}{c}{ $0.519\ (0.127)$ } & $0.153\ (0.067)$ &
$0.175$ \\
&& $\sec$ &\multicolumn{2}{c}{ $0.164\ (0.063)$ } & $0.124\ (0.054)$ &
$0.115$ \\
&& $\thi$ &\multicolumn{2}{c}{ $0.490\ (0.156)$ } & $0.165\ (0.058)$ &
$0.147$ \\
&& $\fou$ &\multicolumn{2}{c}{ $0.206\ (0.072)$ } & $0.078\ (0.036)$ &
$0.071$ \\ [3pt]
& {$\var(Y|X)$} && Zhu--Zhu & Zhu--Dong--Li & \see & ALB
\\
&& $\fir$ & $0.656\ (0.231)$ & $0.283\ (0.126)$ & $0.125\ (0.051)$ &
$0.116$ \\
&& $\thi$ & $0.843\ (0.197)$ & $0.408\ (0.193)$ & $0.116\ (0.055)$ &
$0.113$ \\[3pt]
&{$M(Y|X)$} & & \multicolumn{2}{c}{ AQE } & \see & ALB
\\
&& $\fiv$ & \multicolumn{2}{c}{$0.029\ (0.012)$} & $0.019\ (0.013)$ &
$0.016$ \\
&& $\six$ & \multicolumn{2}{c}{$0.087\ (0.022)$} &$0.049\ (0.019)$ &
$0.043$ \\ [6pt]
{500}& {$E(Y|X)$} & &\multicolumn{2}{c}{ \rmave} & \see & ALB \\
&& $\fir$ &\multicolumn{2}{c}{ $0.100\ (0.030)$ } & $0.081\ (0.021)$ &
$0.083$ \\
&& $\sec$ &\multicolumn{2}{c}{ $0.081\ (0.018)$ } & $0.073\ (0.013)$ &
$0.073$ \\
&& $\thi$ &\multicolumn{2}{c}{ $0.095\ (0.033)$ } & $0.047\ (0.015)$ &
$0.046$ \\
&& $\fou$ &\multicolumn{2}{c}{ $0.079\ (0.015)$ } & $0.047\ (0.010)$ &
$0.045$ \\ [3pt]
&{$\var(Y|X)$} & & Zhu--Zhu & Zhu--Dong--Li & \see & ALB
\\
&& $\fir$ & $0.315\ (0.066)$ & $0.219\ (0.034)$ & $0.109\ (0.020)$ &
$0.104$ \\
&& $\thi$ & $0.236\ (0.035)$ & $0.183\ (0.037)$ & $0.071\ (0.025)$ &
$0.066$ \\ [3pt]
&{$M(Y|X)$} & & \multicolumn{2}{c}{AQE} & \see & ALB
\\
&& $\fiv$ & \multicolumn{2}{c}{$0.017\ (0.005)$} & $0.009\ (0.003)$ &
$0.010$ \\
&& $\six$ & \multicolumn{2}{c}{$0.042\ (0.015)$} & $0.031\ (0.009)$ &
$0.027$ \\
\hline
\end{tabular*}
\end{table}

(c) \textit{Comparison for the central median subspace.}
Let $T(\eta(X, \cdot))$ be the conditional median $M(Y|X)$. We
compare \see\  with the adaptive quantile estimator (AQE) introduced
by \citet{KonXia12}, which can also be used to estimate
the central median subspace.
We use the following models:
\[
\model \ \fiv\dvtx\quad Y = X\lo1\udex2 + X\lo2 \varepsilon,\qquad \model \ \six\dvtx\quad Y = 3 X
\lo1 + X\lo2 + \varepsilon,
\]
where $X \sim N(0, I\lo{10})$ and $X \indep\varepsilon$. For Model
$\fiv$, $\varepsilon$ has a skewed-Laplace distribution with p.d.f.
\[
 f(\varepsilon) = \cases{
(5/4) e\udex{-5 \varepsilon/2}, &\quad $\varepsilon
\geq- (2/5) \log (4/3),$\vspace*{2pt}
\cr
({80}/{27}) e\udex{5 \varepsilon}, &\quad $
\varepsilon< - ({2}/{5}) \log({4}/{3}).$}
\]
In this case,
\[
E(Y|X) = X \lo1 \udex2 + X \lo2 \bigl[1/5 - 2/5 \log(4/3)\bigr],\qquad M(Y|X) = X \lo1
\udex2.
\]
It follows that
\begin{eqnarray*}
\cts E &=&\operatorname{span} \bigl\{ (1, 0, \ldots, 0)\trans, (0, 1, 0, \ldots,
0)\trans\bigr\},\\
 \cts M &=& \operatorname{span} \bigl\{ (1, 0, \ldots, 0)\trans
\bigr\}.
\end{eqnarray*}
For Model $\six$, $\varepsilon\sim t \lo{(3)}$.
Although for this model the central mean subspace coincides with the
central median subspace, due to the heavy-tailed error distribution
the conditional median is preferred to
the conditional mean. Similar models can be found, for example, in
\citet{ZouYua08}.
The results for sample sizes $n = 200, 500$ are presented in Table~\ref{tab2},
in the blocks indicated by $M(Y|X)$.

\def\cms{{{CMS}}}

(d) \textit{Comparison with theoretical lower bound.}
To see how closely the theoretical asymptotic lower bound (\abb{alb})
is approached by \abb{see} for finite samples, we now compute the limit
%
\begin{equation}
\lim\lo{n \to\infty} \sqrt n E \bigl( \| \Pi_ {\mathrm{span}(\hat\beta)
} - \Pi_ {\mathrm{span}(\beta\lo0) } \|
\lo2\bigr), \label{eq:alb}
\end{equation}
where $\hat\beta$ is the semiparametrically efficient estimate. This
is the best we can do to estimate the $T$-central subspace.
The explicit form and the derivation of (\ref{eq:alb}) is given in the supplementary material [\citet{supp}]
(Section \abb{iv}). We present the numerical values
of this limit in the last column (under
the heading \abb{alb}) of Table~\ref{tab2}
for different models and $T$ functionals.

(e) \textit{Conclusions for comparisons in \textup{(a)--(d)}.}
From Table~\ref{tab2}, we see that \see\ achieves
substantially improved accuracy across all models and functionals considered.
Stability of the estimates is also improved as can be seen
from the decrease in standard errors. Our simulation studies (not
presented here)
indicate that the results are not significantly affected by the working
dimension of the central subspace. For example, we repeated
the analysis with $d = 2, 4$ and the patterns of the comparisons are
not significantly altered.

Comparing the results for $n=200$ and $n = 500$, we see that the
proportion of improvement is smaller for the large sample size, as to
be expected.

We see that the actual errors of the \abb{see} computed from
simulations are very close
to the theoretical lower bounds both sample sized $n = 200, 500$ and
the differences become negligible for $n=500$.
Since in the estimator the central subspace is estimated from the
sample and in the lower bounds, the central subspace is treated as known;
the closeness of these errors to their corresponding \abb{alb} also
indicates that the lower bound based on $S\eff\hi*$ in Section~\ref{section:effect} is close to the lower
bound based on $S\eff$. In other words, the effect of estimating the
central subspace on the efficient score is small.

(f) \textit{Comparison under dependent components of $X$}. We now
repeat comparisons in (a) through (d)
using an $X$ with dependent components. Rather than taking $\var(X) =
I \lo{10}$,
we now take
%
\begin{equation}
\label{eq:dependent case} \cov(X \hi i, X \hi j) = 0.5 \hi{| i - j |},\qquad i,j=1, \ldots,10.
\end{equation}
The same covariance matrix was used in \citet{MaZhu12}.
The results parallel to those in Table~\ref{tab2} are presented in Table~\ref{tab3}. We
see that the errors are larger than those for $X$ with independent components,
but the degree by which \abb{see} improves upon the other estimators,
and to which it approaches theoretical asymptotic lower bound, are
similar to those for
the independent-component case.

%
\begin{table}
\caption{Comparison of \see\ with other estimators with correlated
predictors}\label{tab3}
\begin{tabular*}{\textwidth}{@{\extracolsep{\fill}}lccccc@{}}
\hline
{\textbf{Functionals}} & {\textbf{Models}} &
\multicolumn{4}{c}{\textbf{Estimators}} \\
\hline
{$E(Y|X)$} & &\multicolumn{2}{c}{ \rmave} & \see& \abb
{ALB} \\
& $\fir$ &\multicolumn{2}{c}{ $0.520\ (0.155)$ } & $0.164\ (0.066)$ &
$0.168$ \\
& $\sec$ &\multicolumn{2}{c}{ $0.404\ (0.165)$ } & $0.160\ (0.080)$ &
$0.149$\\
& $\thi$ &\multicolumn{2}{c}{ $0.571\ (0.134)$ } & $0.304\ (0.075)$ &
$0.289$ \\
& $\fou$ &\multicolumn{2}{c}{ $0.283\ (0.090)$ } & $0.075\ (0.023)$ &
$0.085$ \\[3pt]
{$\var(Y|X)$} & & Zhu--Zhu & Zhu--Dong--Li & \see & \abb
{ALB} \\
& $\fir$ & $0.539\ (0.174)$ & $0.431\ (0.230)$ & $0.131\ (0.077)$ &
$0.108$ \\
& $\thi$ & $0.617\ (0.222)$ & $0.303\ (0.169)$ & $0.204\ (0.066)$ &
$0.227$ \\[3pt]
{$M(Y|X)$} & & \multicolumn{2}{c}{AQE} & \see& \abb
{ALB} \\
& $\fiv$ & \multicolumn{2}{c}{$0.074\ (0.023)$} & $0.015\ (0.012)$ &
$0.012$\\
& $\six$ & \multicolumn{2}{c}{$0.076\ (0.061)$} & $0.032\ (0.015)$ &
$0.041$ \\
\hline
\end{tabular*}
\end{table}
(g) \textit{Comparison for conditional upper quartile.}
We now apply \see\ to estimating the central upper-quartile subspace
in which the functional of interest is solution to the equation
$P ( Y \le c | X) = 0.75$. We generate $X$ from $N (0, \Sigma)$ with
$\Sigma$ given by~(\ref{eq:dependent case}).
We compare \abb{see} with \abb{aqe} for Models
$\fiv$ and $\six$, and the additional model
\[
\model\ \sev\dvtx\quad Y = 1 + X\lo1 + (1 + 0.4 X\lo2) \varepsilon,
\]
where $\varepsilon\sim N (0, 1)$. In Model $\fiv$, the central
upper-quartile subspace has dimension~2, spanned by $(1, 0, \ldots,
0)\trans$ and
$(0, 1, 0, \ldots, 0)\trans$; in Models $\six$ and~$\sev$, the central
upper-quartile subspaces have dimension $1$ and are spanned by $(3, 1,
0, \ldots, 0)\trans$ and $(1, 0.4 \Phi\inv(0.25), 0,\ldots, 0)\trans$,
respectively, where $\Phi$ is the c.d.f. of the standard normal distribution.
The performance of the estimators is summarized in Table~\ref{tab4}.

We see that \see\ outperforms AQE both in average accuracy and
estimation stability.
It is also interesting to note that, for Model $\six$, the central
median subspace coincides with the central upper-quartile subspace,
and the \abb{see} based on the conditional median (Table~\ref{tab3}) performs
better than the \abb{see} based on the conditional upper quartile
(Table~\ref{tab4}).

\begin{table}
\caption{Comparison of \see\ with AQE at $\tau= 0.75$}\label{tab4}
\begin{tabular*}{\textwidth}{@{\extracolsep{\fill}}lccc@{}}
\hline
{\textbf{Models}} & \textbf{AQE} & \textbf{\see} & \textbf{ALB} \\
\hline
$\fiv$ & {$0.176\ (0.083)$} & $0.043\ (0.017)$ & $0.042$ \\
$\six$ & {$0.160\ (0.037)$} &$0.069\ (0.021)$ & $0.053$ \\
$\sev$ & {$0.245\ (0.117)$} &$0.109\ (0.072)$ & $0.111$ \\
\hline
\end{tabular*}
\end{table}

\section{Application: Age of abalones}\label{section:real data}

\def\cvs{{{CVS}}}

In this section, we evaluate the performance of \see\  in an
application, which is concerned with predicting the age of abalones
using their physical measurements. The data can be found at the
website \url{http://archive.ics.uci.edu/ml/datasets.html}, and consist
of observations from 4177 abalones, of which 1528 are male, 1307 are female
and 1342 are infant. The observations on each subject contain 7
physical measurements and the age of the subject, as measured by the
number of rings in its shell.
We only use the subset of male abalones. The 532th subject in this
subset is an outlier, and is deleted.
Thus, we have a sample of size $1527$ with $7$ predictors and 1
response. For objective evaluation of the estimators,
we further split the data into two subsets: the first 764 subjects are
used as the training set to
estimate the sufficient predictors and the rest 763 subjects are used
as the testing set to plot the derived sufficient predictors versus the
response.

\begin{figure}

\includegraphics{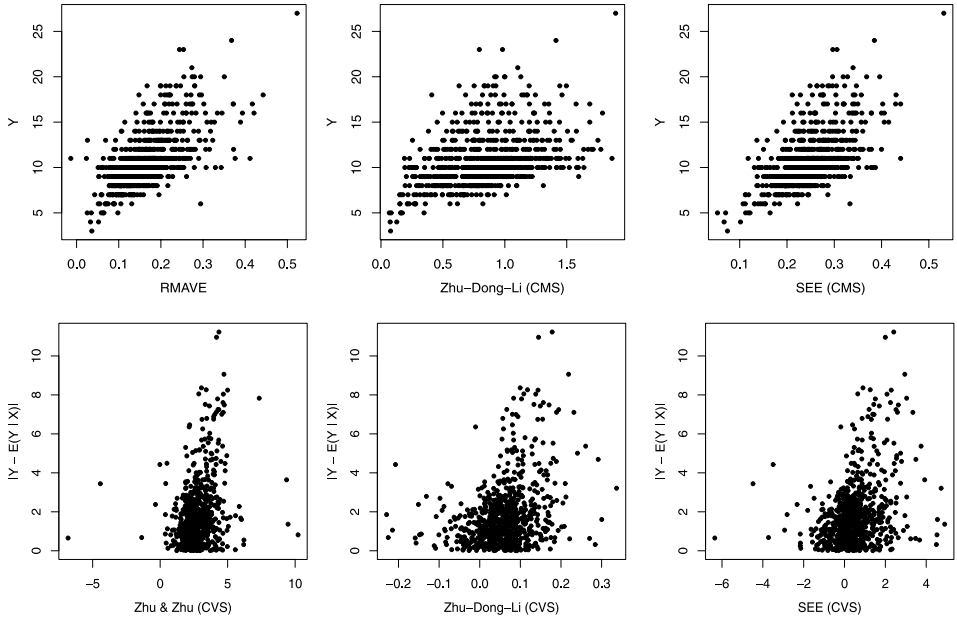}

\caption{Comparison of \see\ with other estimators of \cms\ and \cvs\
for the abalone data.
In the upper panels, the x-axes are the predictors obtained by \rmave,
Zhu--Dong--Li and \see\  estimates for the \cms; the y-axes
are the abalones' ages. In the lower panels,
the x-axes are the predictors derived from
Zhu--Zhu, Zhu--Dong--Li and \see\ estimates of the \cvs; the y-axes are
the estimated absolute residuals.}
\label{figure:abalone}
\end{figure}

We estimate both the central mean subspace (\cms) and the central
variance subspace (\cvs) of this data set.
The \cms\ is estimated by \rmave, \see, and the method implicitly
contained in Zhu, Dong and Li (\citeyear{ZhuDonLi13}).
The \cvs\ is estimated by the methods proposed by \citet{ZhuZhu09}
and Zhu, Dong and Li (\citeyear{ZhuDonLi13}), and the \see.
The results are presented in Figure~\ref{figure:abalone}. The three
upper panels are scatter plots of $Y$ versus the sufficient predictor
in the \cms\ as estimated by \rmave, Zhu--Dong--Li and \see\ in that
order. The three lower panels are the scatter plots of the absolute
residuals $|Y - \hat E (Y|X)|$ versus the sufficient predictor in the
\cvs\ as estimated by Zhu--Zhu, Zhu--Dong--Li and \see.

To give an objective numerical comparison, we use a bootstrapped error
measurement akin to that introduced
by \citet{YeWei03}, which is reasonable because all estimators
involved are consistent.
Since the predictors in the abalone data set
are highly correlated, two estimates of $\beta$ that span substantially
different linear spaces can correspond\vspace*{1pt} to nearly identical $\beta\trans X$.
For this reason, rather than measuring the error in $\hat\beta$, as we
did in the simulations,
here we directly measure the error in $\hat\beta\trans X$.
Specifically, we generate $500$ bootstrap samples, and for each sample
we compute the estimate $\tilde\beta$.
We also compute the full-sample estimate $\hat\beta$. For each
bootstrap sample,
we evaluate the sample correlation between
\[
\bigl\{ \tilde\beta\trans X \lo1, \ldots, \tilde\beta\trans X \lo n\bigr\},
\qquad\bigl\{ \hat\beta
\trans X \lo1, \ldots, \hat\beta\trans X \lo n \bigr\}.
\]
We denote sample correlations for the 500 bootstrap samples as $\rho\lo
1,\ldots, \rho\lo{500}$. We then compute
$1 - \sum\lo{i=1}\udex{500} |\rho\lo i| / 500$ and call it the
bootstrap error of the estimator. The result is summarized in the
Table \ref{tab5}.

\begin{table}
\caption{Bootstrap error of the estimators}\label{tab5}
\begin{tabular*}{\textwidth}{@{\extracolsep{\fill}}lccc@{}}
\hline
{\textbf{Functionals}} & \multicolumn{3}{c@{}}{\textbf{Estimators}} \\
\hline
{$E(Y | X)$} & \rmave& Zhu--Dong--Li & \see\\
& $0.145$ & $0.009$ & $0.003$ \\[3pt]
{$\var(Y | X)$} & Zhu--Zhu & Zhu--Dong--Li & \see\\
& $0.213$ & $0.131$ & $0.105$ \\
\hline
\end{tabular*}\vspace*{-6pt}
\end{table}

We see that \abb{see} is the top performer for estimating both \abb
{cms} and \abb{cvs}, followed by the estimator of Zhu, Dong and Li
(\citeyear{ZhuDonLi13}), and then by
\abb{RMAVE}.
We also observe that the estimation of central variance subspace is in
general less accurate than that of central mean subspace, as has been
observed in many other cases, for example, in Zhu, Dong and Li
(\citeyear{ZhuDonLi13}).

\section{Discussions}\label{section:discussion}

In this paper, we introduce a general paradigm for sufficient dimension
reduction with respect to a conditional
statistical functional, along with semiparametrically efficient
procedures to estimate the sufficient predictors of that functional.
This method is particularly useful when we want to select sufficient predictors
with some specific purposes in mind, such as estimating the conditional
quantiles in a population. This work is a continuation,
synthesis and refinement of previous works on nonparametric mean
regression, nonparametric quantile regression and nonparametric
estimation of heteroscedasticity,
under the unifying framework of \abb{sdr}.
It provides us with tools to explore the detailed structures of the
central subspace, making \abb{sdr} more specific to our goals. Our work
has also substantially broadened the scope of
the semiparametric approach recently introduced to \abb{sdr} by Ma and
Zhu (\citeyear{MaZhu12,MaZhu13N1} and \citeyear{MaZhu13N2}).

In a wide range of simulation studies, the \abb{see} is shown to
outperform several previously proposed estimators for conditional mean,
conditional quantile
and conditional variance. Moreover, the theoretical semiparametric
lower bound is approximately achieved by the actual error based on simulation.
Finally, the algorithm we developed for \abb{see} has a special
advantage over that proposed in \citet{MaZhu13N1}: it does not
rely on any specific parameterization
of the central subspace, which means we do not need to subjectively
assign any element of $\beta$ to be nonzero from the outset.

\section*{Acknowledgments}
The authors would like to thank two referees and an Associate Editor
for their
insightful and constructive comments and suggestions, which led to
significant improvement of
this work.

\begin{supplement}[id=suppA]
\stitle{External appendix to
``On efficient dimension reduction with respect to a
statistical functional of interest''}
\slink[doi]{10.1214/13-AOS1195SUPP} 
\sdatatype{.pdf}
\sfilename{aos1195\_supp.pdf}
\sdescription{The supplementary file provides the proof of
Theorem \ref{theorem:tangent eta}, explicit formulas for the efficient scores
in Section \ref{section:simulation}, the
efficient score when the central subspace
is unknown, and the explicit value of the limit in (\ref{eq:alb}).}
\end{supplement}

%


\printaddresses

\end{document}